\pdfoutput=1
\documentclass[a4paper]{article}
\usepackage[margin=30mm]{geometry}
\usepackage[svgnames]{xcolor}
\usepackage{amsmath, amssymb, amsthm, mathtools, comment}
\usepackage{authblk}
\usepackage[normalem]{ulem}
\usepackage[setpagesize=false]{hyperref}
\hypersetup{
 colorlinks=true,
 linkcolor=FireBrick,
 citecolor=DarkBlue,
}
\usepackage{tikz}
\usetikzlibrary{cd}

\title{Center construction for group-crossed tensor categories}
\author{Mizuki Oikawa\thanks{The author is supported by JSPS KAKENHI Grant Number JP23KJ0540.}}
\affil{Graduate School of Mathematical Sciences \\ The University of Tokyo, Komaba, Tokyo 153-8914, Japan \\ E-mail: \texttt{moikawa@ms.u-tokyo.ac.jp}}
\date{}

\theoremstyle{definition}
\newtheorem{definition}{Definition}[section]
\newtheorem{theorem}[definition]{Theorem}
\newtheorem{lemma}[definition]{Lemma}
\newtheorem{proposition}[definition]{Proposition}
\newtheorem{corollary}[definition]{Corollary}
\newtheorem{remark}[definition]{Remark}
\newtheorem{example}[definition]{Example}
\newtheorem{task}[definition]{Task}

\makeatletter
    
    \@addtoreset{equation}{section}
  \makeatother

\providecommand{\abs}[1]{\left\lvert#1\right\rvert}

\newcommand\cala{\mathcal{A}}
\newcommand\calb{\mathcal{B}}
\newcommand\calc{\mathcal{C}}
\newcommand\cald{\mathcal{D}}
\newcommand\cale{\mathcal{E}}
\newcommand\calm{\mathcal{M}}
\newcommand\caln{\mathcal{N}}
\newcommand\id{\mathrm{id}}
\renewcommand\hom{\operatorname{Hom}}

\DeclareMathOperator{\ad}{Ad}
\DeclareMathOperator{\aut}{Aut}
\DeclareMathOperator{\bimod}{Bimod}
\DeclareMathOperator{\bfmod}{\mathbf{Mod}}
\DeclareMathOperator{\coend}{Coend}

\DeclareMathOperator{\en}{End}
\DeclareMathOperator{\enbar}{\underline{End}}
\DeclareMathOperator{\forg}{Forg}
\DeclareMathOperator{\fpdim}{FPdim}
\DeclareMathOperator{\fun}{Fun}
\DeclareMathOperator{\hombar}{\underline{Hom}}
\DeclareMathOperator{\homog}{Homog}
\DeclareMathOperator{\irr}{Irr}
\DeclareMathOperator{\obj}{Obj}
\DeclareMathOperator{\rep}{Rep}
\DeclareMathOperator{\rmz}{Z}
\DeclareMathOperator{\st}{st}
\DeclareMathOperator{\tr}{Tr}

\renewcommand\mod{\operatorname{Mod}}

\begin{document}
\tikzset{>=stealth}
\tikzset{cross/.style={preaction={-,draw=white,line width=6pt}}}
\tikzset{block/.style={draw, fill=white, rectangle, minimum width=0.5}}
\maketitle

\begin{abstract}
	We define the notion of the $(G,\Gamma)$-crossed center of a $(G,\Gamma)$-crossed tensor category in the sense of Natale. We show that the $(G, \Gamma)$-crossed center is a $(G \bowtie \Gamma, G \times \Gamma)$-braided tensor category. This construction generalizes the graded center construction for graded tensor categories and the equivariant center construction for tensor categories with group actions. 
\end{abstract}

\tableofcontents

\section{Introduction}

Groups have been extensively studied as fundamental symmetries in various mathematical and physical contexts. However, it is sometimes important to consider symmetries beyond groups, which can provide a more general framework for understanding various structures. One such example of more general symmetries is the concept of a \emph{tensor category}. A tensor category is a category equipped with a product structure, allowing the composition of objects in a manner that generalizes the concept of group multiplication. An important example of tensor categories is the representation category of conformal field theory \cite{MR297259}\cite{MR1016869}. The product structure in these tensor categories is commutative, leading to what is known as a \emph{braided tensor category}. Braided tensor categories are crucial for understanding the correlation functions in conformal field theory \cite{MR1940282}. Moreover, from braided tensor categories, one can construct the Reshetikhin--Turaev topological quantum field theory \cite{MR1091619}, which highlights the importance of these structures from a topological perspective as well.

When a conformal field theory admits a group action, it is possible to consider the equivariant version of its representation category \cite{MR2183964}. Although this category is not braided in the usual sense, it possesses a structure called a \emph{$G$-braiding} (in the sense of Turaev \cite{turaev2000homotopy}), which can be regarded as an equivariant version of braiding. A $G$-braided tensor category is a tensor category equipped with both a compatible $G$-action and $G$-grading. These structures modify the commutativity of the product in a specific way. The categorical fixed points under the $G$-action in a $G$-braided tensor category yield a braided tensor category in the usual sense. In the case of a finite group action on a rational conformal field theory, the categorical fixed points of the category of twisted representations correspond to the ordinary representation category of the fixed points in the conformal field theory \cite{MR2183964}. This fact allows us to describe the representation category of the finite group fixed points in the original conformal field theory. In general, while defining group fixed points is straightforward, providing an explicit description is often challenging, making this result particularly useful.

Furthermore, through the equivariant version of the Reshetikhin--Turaev construction, one can construct a \emph{homotopy quantum field theory} from a $G$-braided tensor category \cite{turaev2000homotopy} \cite{MR3195545}, which serves as an invariant of manifolds with principal $G$-bundles. This offers intriguing insights from a topological perspective. Thus, a $G$-braided tensor category is interesting and significant for two main reasons: it provides a fascinating invariant in its own right, and it yields a usual braided tensor category through categorical fixed points, thereby offering a more refined description of that braided tensor category.

Through the Drinfeld center \cite{MR1107651}\cite{MR1151906}, which serves as the categorical version of the center of a unital noncommutative ring, it is possible to construct a braided tensor category from a tensor category that is not necessarily braided. There are at least two generalizations of the Drinfeld center construction that produce $G$-braided tensor categories. One is the graded center construction for graded tensor categories \cite{MR2587410}\cite{MR3053524}, and another is the equivariant center construction for tensor categories with group actions \cite{MR3936135}\cite{sheikh}. These constructions are important for providing examples of $G$-braided tensor categories; however, they utilize only one of the two structures: either the $G$-action or the $G$-grading. Therefore, it is natural to consider constructions for tensor categories where these structures are compatible, as in the case of $G$-braided tensor categories.

In this article, we define the notion of the \emph{$(G,\Gamma)$-crossed center} of a $(G,\Gamma)$-crossed pivotal tensor category in the sense of Natale \cite{MR3550925}. This construction reproduces the graded center construction when the group $G$ acting on the tensor category is trivial. Moreover, it reproduces the equivariant center construction when the group $\Gamma$ that grades the tensor category is trivial. 

We show that the $(G,\Gamma)$-crossed center of a $(G,\Gamma)$-crossed pivotal tensor category is a \emph{$(G \bowtie \Gamma, G \times \Gamma)$-braided} tensor category, where $G \bowtie \Gamma$ denotes the \emph{Zappa--Sz\'{e}p product} \cite{MR0019090} of $G$ and $\Gamma$. Indeed, our construction produces a group-crossed tensor category in the sense of Natale even when we start from one in the sense of Turaev.

This article is structured as follows. In Section \ref{section_preliminaries}, we review some preliminaries on group-crossed and group-braided monoidal categories in the sense of Natale. In Section \ref{subsection_coherence}, we show the coherence theorem for $(G,\Gamma)$-crossed monoidal categories. In Section \ref{section_center}, we define the notion of the \emph{$(G,\Gamma)$-crossed center} of a $(G,\Gamma)$-crossed pivotal tensor category. 

\section{Preliminaries}
\label{section_preliminaries}

In this section, we recall the notion of a $(G,\Gamma)$-crossed monoidal category in the sense of Natale \cite{MR3550925}. In short, it is a monoidal category with a compatible group action and a grading in some sense. We also recall the counterpart of a braiding on such a category.

Our terminologies on tensor category theory follow those in \cite{egno}. In particular, a \emph{multitensor category} over a field $k$ is a locally finite $k$-linear abelian rigid monoidal category with a bilinear monoidal product as in \cite[Definition 4.1.1]{egno}. A \emph{tensor category} is a multitensor category such that the endomorphism space of the unit is one-dimensional. Here, an \emph{abelian monoidal category} is a monoidal category that is an abelian category with a biexact monoidal product. A \emph{linear monoidal category} is a monoidal category that is a linear category (i.e. a category with linear hom spaces whose composition of morphisms are bilinear) with a bilinear monoidal product.

First, we recall the notion of a \emph{grading} on an additive category. Let $\{ \cala_i \}_{i \in I}$ be a family of additive categories. Then, its \emph{direct sum} $\calc = \bigoplus_{i \in I} \cala_i$ is an additive category with a family of additive functors $\{ I_i : \cala_i \to \cala \}_{i \in I}$ which has the universality that for any additive category $\calb$ and a family of additive functors $\{ F_i: \cala_i \to \calb \}_{i \in I}$ there exists an additive functor $F: \cala \to \calb$ with a family of natural isomorphisms $\sigma^F = \{ \sigma^{F,i}: F I_i \cong F_i \}_{i \in I}$ such that for any other such pair $(\tilde{F},\tilde{\sigma}^F)$ and a family of natural transformations $\{ \tau^i : F_i \to \tilde{F}_i \}_{i \in I}$ we have a unique natural transformation $\tau : F \to \tilde{F}$ with $\tau^i \sigma^{F,i} = \sigma^{\tilde{F},i} (\tau \ast I_i)$ for any $i \in I$, where $\ast$ denotes the horizontal composition of natural transformations. Such a category $\cala$ always exists, see e.g. \cite[Section 1.3]{egno}. The datum $(\{ \cala_i \}_{i \in I}, \{ I_i \}_{i \in I})$ is called an \emph{$I$-grading} on $\cala$ and an additive category with an $I$-grading is called an \emph{$I$-graded category}.
 
We can indeed take $\sigma^{F, i}$'s to be identities. In this case, we refer to the functor $F$ as the \emph{extension} of the functors $\{ F_i \}_{i \in I}$. Also, we refer to the natural transformation $\tau$ obtained by universality as the \emph{extension} of the natural transformations $\{ \tau^i \}_{i \in I}$. Thus, we often define functors on $\calc$ and natural transformations between them by defining them on $\cala_i$'s. 

More explicitly, the notion of an $I$-grading on an additive category is equivalent to a family of additive subcategories $\{ \cala_i \}_{i \in I}$ such that every object is a direct sum of objects of $\cala_i$'s and $\hom_\cala(\lambda, \mu) = \{ 0 \}$ for $\lambda \in \obj(\cala_i)$ and $\mu \in \obj(\cala_j)$ with $i \neq j$. A decomposition of an object into objects of $\cala_i$'s is called a \emph{homogeneous decomposition}, which is essentially unique by \cite[Lemma 2.10]{oikawa2023frobenius}.

We write $\partial_\cala \lambda = i$ or simply $\partial \lambda =i$ if $\lambda \in \obj(\cala_i)$. Since $\obj(\cala_i) \cap \obj(\cala_j)$ consists of the zero objects of $\cala$ for $i \neq j$, the assignment $\partial$ is single-valued on nonzero objects. We put $\homog(\cala) \coloneqq \bigcup_{i \in I} \obj(\cala_i)$.

We are interested in the case where $\cala$ is a monoidal category, $I$ is a group and they are compatible in an appropriate sense. 

\begin{definition}
	\label{definition_graded_monoidal}
	Let $\Gamma$ be a group. A \emph{$\Gamma$-graded monoidal category} is an additive monoidal category $\calc$ (i.e. a monoidal category that is an additive category with a biadditive monoidal product) that is $\Gamma$-graded and $\obj(\calc_s) \otimes \obj(\calc_t) \subset \obj(\calc_{st})$ for all $s, t \in \Gamma$. 

	A \emph{$\Gamma$-graded monoidal functor} $F: \calc \to \cald$ between $\Gamma$-graded monoidal category is an additive monoidal functor such that $F (\obj(\calc_s)) \subset \obj(\cald_s)$ for all $s \in \Gamma$. 
\end{definition}

Next, we recall the notion of an \emph{action} of a group on a (not necessarily monoidal) category.

\begin{definition}
	An \emph{action} of a group $G$ on a category $\cale$ is a monoidal functor from the monoidal category $\underline{G}$ of elements of $G$ with only identity morphisms to the monoidal category of endofunctors on $\cale$. 

	When $\cale$ is additive (resp. linear), we require that an action of a group values in the category of additive (resp. linear) endofunctors on $\cale$. 
\end{definition}

More explicitly, an action $\gamma$ of $G$ on $\cale$ consists of an assignment $g \mapsto \gamma(g)$ of endofunctors and natural isomorphisms $\chi^\gamma_{g,h}: \gamma(g) \gamma(h) \cong \gamma(g h)$ for $g, h \in G$ and $\iota^\gamma: \id_\cale \cong \gamma(e)$ with certain coherence relations, see e.g. \cite[Section 2.4]{egno}. By the coherence theorem for monoidal functors (see e.g. \cite[Subsection 2.3.3]{MR3076451}), a natural isomorphism obtained by vertical and horizontal compositions of $\chi^\gamma_{g,h}$'s and $\iota^\gamma$ is canonical, see also \cite[Proposition 4.1]{MR3671186}. Therefore, we often suppress such a natural transformation and its components.

Here, we introduce some graphical notations in the bicategory of categories. Let the natural isomorphisms $\chi^{\gamma}_{g,h}$ and $\iota^\gamma$, respectively, be graphically denoted by a fork from $\gamma(g) \gamma(h)$ to $\gamma(g h)$ and a half line with a circle from $\id$ to $\gamma(e)$ like the graphical notation for an algebra (see e.g. \cite[Subsection 2.3]{}). Then, the coherence relations for $\chi^\gamma$ and $\iota^\gamma$ graphically correspond to the axioms for an algebra. 

For $\lambda \in \obj(\cale)$ and $g \in G$, we let ${}^g \lambda$ denote $\gamma(g)(\lambda)$. We also use a similar notation for morphisms. 

In order to formulate compatibility between a group action and a grading, we need the following notion, see e.g. \cite[Definition IX.1.1]{MR1321145}.

\begin{definition}
	\label{definition_matched_pair}
	A \emph{matched pair of groups} is a tuple $(G,\Gamma) = (G,\Gamma, \rhd_1, \rhd_2)$ of groups $G$ and $\Gamma$, an action $\rhd_1: G \times \Gamma \to \Gamma$ of $G$ on a set $\Gamma$ and an action $\rhd_2: \Gamma \times G \to G$ of $\Gamma$ on a set $G$ such that $g \rhd_1 e = e$, $s \rhd_2 e = e$ and
	\begin{align*}
		g \rhd_1 (st) &= ((t \rhd_2 g) \rhd_1 s)(g \rhd_1 t) \\
		s \rhd_2 (gh) &= ((h \rhd_1 s) \rhd_2 g)(s \rhd_2 h)
	\end{align*} 
	for all $g,h \in G$ and $s,t \in \Gamma$.
\end{definition}

The notion of a matched pair of groups is equivalent to that of an \emph{exact factorization} of a group. Namely, an \emph{exact factorization} of a group $H$ is a pair $(G,\Gamma)$ of subgroups with $H = G \Gamma$ and $G \cap \Gamma = \{e\}$. Given an exact factorization $(G,\Gamma)$ of a group $H$, we can obtain a matched pair $(G,\Gamma,\rhd_1,\rhd_2)$ from the relation $s g = (s \rhd_2 g^{-1})^{-1} (g^{-1} \rhd_1 s)$ for $g \in G$ and $s \in \Gamma$. Conversely, given a matched pair $(G,\Gamma,\rhd_1,\rhd_2)$, this relation defines the \emph{Zappa--Sz\'{e}p product} $H \coloneqq G \bowtie \Gamma$ of $G$ and $\Gamma$ \cite{MR0019090} (see also \cite[Proposition IX.1.2]{MR1321145}), which has an exact factorization $(G,\Gamma)$.

\begin{remark}
	\label{remark_right_action}
	In the usual definition of a matched pair (e.g. \cite[Definition IX.1.1]{MR1321145}), $\rhd_1$ is a right action. Our definition is obtained by regarding it as a left action by the inverse map $G^{\mathrm{op}} \cong G$.
\end{remark}

We are interested in the following notion, which was introduced in \cite[Definition 4.1]{MR3550925}.

\begin{definition}
	\label{definition_crossed_category}
	Let $(G, \Gamma)$ be a matched pair of groups. A \emph{$(G, \Gamma)$-crossed monoidal category (in the sense of Natale)} is a $\Gamma$-graded additive monoidal category $\calc = \bigoplus_{s \in \Gamma} \calc_s$ with an action $\gamma^\calc$ of $G$ on a category $\calc$, a family of isomorphisms $J^{\gamma^\calc} = \{ J^{\gamma^\calc(g)}_{\lambda, \mu}: {}^{\partial \mu \rhd_2 g}\lambda {}^g \mu \cong {}^g (\lambda \mu) \}_{\lambda \in \obj(\calc), \mu \in \homog(\calc), g \in G}$ that is natural in $\lambda$ and $\mu$ and a family of isomorphisms $\varphi^{\gamma^\calc} = \{ \varphi^{\gamma^\calc(g)}: \mathbf{1}_\calc \cong {}^g \mathbf{1}_\calc \}_{g \in G}$ with the following axioms: 
	\begin{enumerate}
		\item For all $s \in \Gamma$ and $g \in G$, $\gamma^\calc(g) (\obj(\calc_s)) \subset \obj(\calc_{g \rhd_1 s})$.
		\item For all $g \in G$, $\lambda \in \obj(\calc)$ and $\mu, \nu \in \homog(\calc)$,
		\begin{align*}
			{}^g a^\calc_{\lambda, \mu, \nu} J^{\gamma^\calc(g)}_{\lambda \mu, \nu} (J^{\gamma^\calc(\partial \nu \rhd_2 g)}_{\lambda, \mu} \otimes \id_{{}^g\nu}) &= J^{\gamma^\calc(g)}_{\lambda, \mu \nu} (\id_{{}^{\partial (\mu \nu) \rhd_2 g}\lambda} \otimes J^{\gamma^\calc(g)}_{\mu, \nu}) a^\calc_{{}^{\partial \mu \partial \nu \rhd_2 g} \lambda, {}^{\partial \nu \rhd_2 g} \mu, {}^g \nu} \\
			J^{\gamma^\calc(g)}_{\mathbf{1}_\calc,\mu} (\varphi^{\gamma^\calc(\partial \mu \rhd_2 g)} \otimes \id_{{}^g \mu}) &= \id_{{}^g \mu} = J^{\gamma^\calc(g)}_{\mu, \mathbf{1}_\calc} (\id_{{}^g \mu} \otimes \varphi^{\gamma^\calc(g)}).
		\end{align*}
		\item For all $g, h \in G$, $\lambda \in \obj(\calc)$ and $\mu \in \homog(\calc)$,
		\begin{align*}
			 (\chi^{\gamma^\calc}_{g, h})_{\lambda \mu} {}^g J^{\gamma^\calc(h)}_{\lambda, \mu} J^{\gamma^\calc(g)}_{{}^{\partial \mu \rhd_2 h} \lambda, {}^h \mu} &= J^{\gamma^\calc(gh)}_{\lambda, \mu} ((\chi^{\gamma^\calc}_{(h \rhd_1 \partial \mu) \rhd_2 g, \partial \mu \rhd_2 h})_\lambda \otimes (\chi^{\gamma^\calc}_{g, h})_\mu) \\
			 (\chi^{\gamma^\calc}_{g,h})_{\mathbf{1}_\calc} {}^g \varphi^{\gamma^\calc(h)} \varphi^{\gamma^\calc(g)} &= \varphi^{\gamma^\calc(gh)} \\
			\iota^{\gamma^\calc}_{\lambda \mu} &= J^{\gamma^\calc(e)}_{\lambda,\mu} (\iota^{\gamma^\calc}_\lambda \otimes \iota^{\gamma^\calc}_{\mu}) \\
			\iota^{\gamma^\calc}_{\mathbf{1}_\calc} &= \varphi^{\gamma^\calc(e)}. 
		\end{align*}
	\end{enumerate} 
\end{definition}

\begin{remark}
	\label{remark_crossed_generalization}
	One can see that axioms 2 and 3 give a generalization of $\gamma^\calc$ being an \emph{action} of $G$ on the monoidal category $\calc$ (i.e. an action of $G$ on the category $\calc$ that values in the category of monoidal endofunctors, see e.g. \cite[Definition 2.7.1(ii)]{egno}). Indeed, when $\rhd_2$ is trivial, axiom 2 becomes the monoidality of $\gamma^\calc(g)$, and axiom 3 becomes the monoidality of $\chi^{\gamma^\calc}_{g,h}$ and $\iota^\gamma$.

	In particular, if $\Gamma$ is trivial, then the notion of a $(G, \Gamma)$-crossed monoidal category is equivalent to that of an additive monoidal category with an action of $G$.
	
	On the other hand, if $G$ is trivial, then by definition the notion of a $(G, \Gamma)$-crossed monoidal category is equivalent to that of a $\Gamma$-graded monoidal category.
\end{remark}

For $\lambda, \mu \in \obj(\calc)$ and $g \in G$, we let ${}^g \lambda \mu$ denote $\gamma^\calc(g)(\lambda) \otimes \mu$. We also use a similar notation for morphisms.

We introduce some graphical notations. For every $g \in G$ and $s \in \Gamma$, we can define a functor $\gamma_s^\calc(g): \calc \times \calc_s \to \calc \times \calc_{g \rhd_1 s}$ of by putting $\gamma_s^\calc(g)(\lambda,\mu) \coloneqq ({}^{s \rhd_2 g} \lambda ,{}^g \mu)$ for $\lambda \in \obj(\calc)$ and $\mu \in \obj(\calc_s)$. Then, the isomorphisms $J^{\gamma^\calc(g)}_{\lambda ,\mu}$ form a natural isomorphism $J^{\gamma^\calc(g)}: \otimes \circ \gamma^{\calc}_s(g) \cong \gamma^{\calc}_s(g) \circ \otimes$. We let $J^{\gamma^\calc(g)}$ be graphically denoted by a crossing from $\otimes g$ to $g \otimes$. Then, the first relation of axiom 3 in Definition \ref{definition_crossed_category} is graphically represented by Figure \ref{figure_axiom3_crossed_category}.

\begin{figure}[htb]
	\centering 
	\begin{tikzpicture}
		\draw[->] (1.25,-0.5) -- (1.25,-0.75) -- (0,-2);
		\draw[->,cross] (0,0) -- (0,-1) -- (1,-2);
		\draw (0.75,0) arc (180:360:0.5);
		\node at (0.75,0.25){$g$};
		\node at (1.75,0.25){$h$};
		\node at (0,0.25){$\otimes$};
		\node at (1,-2.25){$\otimes$};
		\node at (0,-2.25){$g h$};
		\node at (2.5,-1){$=$};
		\begin{scope}[shift={(3.5,0)}]
			\draw[->] (0.5,-1.75) -- (0.5,-2);
			\draw (0.75,0) -- (0,-0.75) -- (0,-1.25) arc (180:360:0.5) -- (1.75,-0.5) -- (1.75,0);
			\draw[->,cross] (0,0) -- (1.75,-1.75) -- (1.75,-2);
			\node at (0.75,0.25){$g$};
			\node at (1.75,0.25){$h$};
			\node at (0,0.25){$\otimes$};
			\node at (1.75,-2.25){$\otimes$};
			\node at (0.5,-2.25){$g h$};
		\end{scope}
	\end{tikzpicture}
	\caption{Axiom 3 for $J^{\gamma^\calc}$}
	\label{figure_axiom3_crossed_category}
\end{figure}

We can also consider a generalization of a braiding, see \cite[Definition 7.1]{MR3550925}.

\begin{definition}
	\label{definition_group_braiding}
	Let $(G, \Gamma) = (G, \Gamma, \rhd_1, \rhd_2)$ be a matched pair of groups. A \emph{braiding} on $(G,\Gamma)$ is the pair $(\phi,\psi)$ of group homomorphisms $\phi, \psi: \Gamma \to G$ such that
	\begin{align*}
		(\phi(s) \rhd_1 t) s &= (\psi(t) \rhd_1 s)t \\
		(s \rhd_2 g) \phi(s) &= \phi(g \rhd_1 s) g \\
		(s \rhd_2 g) \psi(s) &= \psi(g \rhd_1 s) g \\
		s \rhd_2 \phi(t) &= \phi(\psi(s) \rhd_1 t) \\
		s \rhd_2 \psi(t) &= \psi(\phi(s) \rhd_1 t)     
	\end{align*}
	for all $g \in G$ and $s,t \in \Gamma$. A matched pair of groups with a braiding is called a \emph{braided matched pair} of groups.
\end{definition}

\begin{definition}
	\label{definition_crossed_braiding}
	Let $(G,\Gamma) \coloneqq (G,\Gamma,\phi,\psi)$ be a braided matched pair of groups. A \emph{$(G,\Gamma)$-braided monoidal category (in the sense of Natale)} is the pair $(\calc,b^\calc)$ of a $(G,\Gamma)$-crossed monoidal category $\calc$ and a natural isomorphism $b^\calc = \{ b^\calc_{\lambda, \mu}: {}^{\phi(\partial \mu)} \lambda \mu \cong {}^{\psi(\partial \lambda)} \mu \lambda \}_{\lambda, \mu \in \homog(\calc)}$ such that
	\begin{align*}
		&{}^g b_{\lambda, \mu}^\calc J^{\gamma^\calc(g)}_{{}^{\phi(\partial \mu)} \lambda, \mu} ((\chi^{\gamma^\calc}_{\partial \mu \rhd_2 g,\phi(\partial \mu)})^{-1}_{\lambda} (\chi^{\gamma^\calc}_{\phi(g \rhd_1 \partial \mu),g})_\lambda \otimes \id_{{}^g \mu}) \\
		&= J^{\gamma^\calc(g)}_{{}^{\psi(\partial \lambda)} \mu, \lambda} ((\chi^{\gamma^\calc}_{\partial \lambda \rhd_2 g,\psi(\partial \lambda)})_{\mu}^{-1} (\chi^{\gamma^\calc}_{\psi(g \rhd_1 \partial \lambda),g})_{\mu} \otimes \id_{{}^g \lambda}) b_{{}^g \lambda, {}^g \mu}^\calc \\
		&b^\calc_{\lambda \mu, \nu} (J^{\gamma^\calc(\phi(\partial \nu))}_{\lambda, \mu} \otimes \id_\nu) = ((\chi^{\gamma^\calc}_{\psi(\partial \lambda),\psi(\partial \mu)})_{\nu} \otimes \id_{\lambda \mu}) (b^\calc_{\lambda, {}^{\psi(\partial \mu)} \nu} \otimes \id_{\mu}) (\id_{{}^{\partial \mu \rhd_2 \phi (\partial \nu)} \lambda} \otimes b^\calc_{\mu, \nu}) \\
		&b^\calc_{\lambda, \mu \nu} ((\chi^{\gamma^\calc}_{\phi(\partial \mu), \phi(\partial \nu)})_\lambda \otimes \id_{\mu \nu}) = (J^{\gamma^\calc(\psi(\partial \lambda))}_{\mu, \nu} \otimes \id_\lambda)(\id_{{}^{\psi(\phi(\partial \nu) \rhd_1 \partial \lambda)} \mu} \otimes b^\calc_{\lambda, \nu}) (b^\calc_{{}^{\phi(\nu)} \lambda, \mu} \otimes \id_{\nu})
	\end{align*}
	for all $\lambda, \mu, \nu \in \homog(\calc)$ and $g \in G$.
\end{definition}

\begin{remark}
	\begin{enumerate}
		\item The components of a $(G,\Gamma)$-braiding $b^\calc$ are grading preserving by the first axiom for $(\phi, \psi)$ in Definition \ref{definition_group_braiding}. The first (resp. second, resp. third) axiom for $b^\calc$ makes sense by the second and third (resp. fourth, resp. fifth) axioms for $(\phi, \psi)$.
		\item Our definition of $(G,\Gamma)$-crossed and $(G,\Gamma)$-braided monoidal categories coincides with that of \cite{MR3550925} if we take $\rhd_1$ to be a right action, see Remark \ref{remark_right_action}.
	\end{enumerate}
\end{remark}

Indeed, the notion of a $(G,\Gamma)$-crossed monoidal category reproduces as a special case that of a \emph{$G$-crossed monoidal category}, which was introduced by Turaev \cite[Section 2.1]{turaev2000homotopy}, see \cite[Section 8.1]{MR3550925}.

\begin{definition}
	\label{definition_turaev}
	Let $G$ be a group. A \emph{$G$-crossed monoidal category (in the sense of Turaev)} is a $(G, G)$-crossed monoidal category, where $(G, G)$ denotes the matched pair given by the adjoint action $\rhd_1$ and the trivial action $\rhd_2$. 
	A \emph{$G$-braided monoidal category (in the sense of Turaev)} is a $(G, G)$-braided monoidal category, where $(G, G)$ denotes the braided matched pair given by the trivial homomorphism $\phi$ and the identity $\psi$.
\end{definition}

Finally, we give some known examples.

\begin{example}
	For a matched pair $(G,\Gamma)$ of groups, the tensor category $\mathrm{Vec}_\Gamma$ of finite dimensional $\Gamma$-graded vector spaces can trivially be regarded as a $(G,\Gamma)$-crossed tensor category. Namely, for $g \in G$ and a vector space $V$ with grading $s \in \Gamma$, define $\gamma(g)(V)$ to be $V$ with grading $g \rhd_1 s$. Then, $\gamma$ and identity morphisms form a $(G,\Gamma)$-structure on $\mathrm{Vec}_\Gamma$.
\end{example}

Before introducing other examples, let us recall the notion of a \emph{relative center} of a monoidal category \cite[Definition 3.2]{MR1151906}. Let $\calc$ be a monoidal category and let $\cald$ be its monoidal subcategory. Then, the \emph{relative center} $\rmz_{\cald}(\calc)$ is defined as follows. An object of $\rmz_{\cald}(\calc)$ is the pair $(\lambda, h^\lambda)$ of $\lambda \in \obj(\calc)$ and a family $h^\lambda = \{ h^\lambda_\mu :\lambda \mu \cong \mu \lambda \}_{\mu \in \obj(\cald)}$ of isomorphisms that is natural in $\mu$ such that $h^{\lambda}_{\mu_1 \mu_2} = (\id_{\mu_1} \otimes h^\lambda_{\mu_2}) (h^\lambda_{\mu_1} \otimes \id_{\mu_2})$ for all $\mu_1, \mu_2 \in \obj(\cald)$. The monoidal product is given by $h^{\lambda_1 \lambda_2}_\mu \coloneqq (h^{\lambda_1}_\mu \otimes \id_{\lambda_2}) (\id_{\lambda_1} \otimes h^{\lambda_2}_\mu)$ for $(\lambda_1, h^{\lambda_1}), (\lambda_2, h^{\lambda_2}) \in \obj(\rmz_\cald(\calc))$ and $\mu \in \obj(\cald)$. The unit is given by $h^{\mathbf{1}_\calc} \coloneqq \id$. 

We also recall that a \emph{pivotal monoidal category} is a rigid monoidal category $\calc$ with a monoidal isomorphism $\delta^\calc: \id \cong {}^{\vee \vee}$, where ${}^\vee$ denotes the left dual, see e.g. \cite[Definition 4.7.7]{egno}. 

\begin{example}
	\label{example_graded_center}
	Let $\calc$ be a $\Gamma$-graded tensor category. If $\calc$ and $\Gamma$ are finite \cite[Theorem 3.3]{MR2587410} or $\calc$ is (not necessarily finite but) pivotal \cite[Theorem 4.1]{MR3053524}, then the \emph{$\Gamma$-graded center} $\rmz_\Gamma(\calc) \coloneqq \rmz_{\calc_e}(\calc)$ of $\calc$ is a $\Gamma$-braided tensor category in the sense of Turaev.
\end{example}

\begin{example}
	\label{example_graded_center}
	Let $\calc$ be an additive monoidal category with an action $\gamma^\calc$ of a group $G$. Note that the assignment $g \mapsto \gamma^\calc(g^{-1})$ defines an action of $G^{\mathrm{op}}$ on $\calc$. Then, the \emph{$G$-equivariant center} $\rmz^G(\calc) \coloneqq \rmz_{\calc}(G^{\mathrm{op}} \ltimes \calc)^{\mathrm{rev}}$ of $\calc$ is a $G$-braided monoidal category in the sense of Turaev \cite[Theorem 6.2]{MR3936135}\cite[Proposition 4.4.6]{sheikh}.
\end{example}

Here, we do not give the details of these two constructions because we give their generalization in Section \ref{section_center}.

\section{Coherence theorems for crossed monoidal categories}
\label{subsection_coherence}

It is known that a group action on a monoidal category can be \emph{strictified} in an appropriate sense \cite[Theorem 4.3]{MR3671186} (see also \cite[Theorem 3.1]{MR3936135}). In this section, we generalize this fact in the setting of group-crossed monoidal categories in the sense of Natale (Definition \ref{definition_crossed_category}). As a consequence, we can see that isomorphisms obtained from the crossed structure are indeed canonical (Corollary \ref{corollary_coherence}), and therefore we may suppress them. 

First, we define the appropriate notion of an \emph{equivalence} between group-crossed monoidal categories. For a monoidal functor $F: \calc \to \cald$ between monoidal categories, let $J^F$ and $\varphi^F$ denote the monoidal structure of $F$. Namely, for a monoidal functor $F$, we have the isomorphisms $J^F_{\lambda, \mu}: F(\lambda) F(\mu) \cong F(\lambda \mu)$ for $\lambda, \mu \in \obj(\calc)$ and $\varphi^F: \mathbf{1}_\cald \cong F(\mathbf{1}_\calc)$, see e.g. \cite[Definition 2.4.5]{egno}.

\begin{definition}
	\label{definition_crossed_functor}
	Let $\calc$ and $\cald$ be $(G, \Gamma)$-crossed monoidal categories for a matched pair $(G,\Gamma)$ of groups. A \emph{$(G, \Gamma)$-crossed monoidal functor} from $\calc$ to $\cald$ is the pair $F = (F,\eta^F)$ of a $\Gamma$-graded monoidal functor $F: \calc \to \cald$ and a family of natural isomorphisms $\eta^F = \{ \eta^F_g: F \gamma^\calc(g) \cong \gamma^\cald (g) F \}_{g \in G}$ such that
	\begin{enumerate}
		\item $\eta^F_{g h} (\id_F \ast \chi^{\gamma^\calc}_{g,h}) = (\chi^{\gamma^\cald}_{g,h} \ast \id_F) (\id_{\gamma^\cald(g)} \ast \eta^F_h)(\eta^F_g \ast \id_{\gamma^\calc(h)})$ for all $g,h \in G$.
		\item $(\eta^F_g)_{\lambda \mu} F(J^{\gamma^\calc(g)}_{\lambda,\mu}) J^F_{{}^{\partial \mu \rhd_2 g} \lambda, {}^g \mu} = {}^g J^F_{\lambda, \mu} J^{\gamma^\cald(g)}_{F(\lambda),F(\mu)} ((\eta^F_{\partial \mu \rhd_2 g})_\lambda \otimes (\eta^F_g)_\mu)$ for all $\lambda \in \obj(\calc), \mu \in \homog(\calc)$ and $g \in G$.
	\end{enumerate}

	Let $F, F' : \calc \to \cald$ be $(G, \Gamma)$-crossed monoidal functors. A \emph{$(G,\Gamma)$-crossed natural transformation} from $F$ to $F'$ is a monoidal natural transformation $\tau: F \to F'$ such that $(\id_{\gamma^\cald(g)} \ast \tau) \eta_g^F = \eta^{F'}_g (\tau \ast \id_{\gamma^{\calc}(g)})$ for all $g \in G$. 
\end{definition}

\begin{remark}
	When $G$ is trivial, the notion of a $(G,\Gamma)$-crossed functor is equivalent to that of a $\Gamma$-graded functor (Definition \ref{definition_graded_monoidal}) by taking $\eta^F$ to be the identity transformation. When $\Gamma$ is trivial, the notion of a $(G,\Gamma)$-crossed functor is equivalent to that of an additive \emph{monoidal $G$-functor}, see e.g. \cite[Definition 2.8]{oikawa2023frobenius}. When $\calc$ and $\cald$ are $G$-crossed monoidal categories in the sense of Turaev (Definition \ref{definition_turaev}), the notion of a $(G,\Gamma)$-crossed functor is equivalent to that of a \emph{$G$-crossed monoidal functor}, see \cite[Definition 2.11]{oikawa2023frobenius}.
\end{remark}

\begin{proposition}
	\label{proposition_crossed_bicategory}
	Let $(G, \Gamma)$ be a matched pair of groups. Then, $(G, \Gamma)$-crossed monoidal categories, $(G, \Gamma)$-crossed monoidal functors and $(G, \Gamma)$-crossed monoidal natural transformations form a 2-category (i.e. strict bicategory).

	\begin{proof}
		It is easy to see that for $(G,\Gamma)$-crossed monoidal categories $\calc$ and $\cald$, $(G,\Gamma)$-crossed monoidal functors from $\calc$ to $\cald$ and $(G,\Gamma)$-crossed monoidal natural transformations between them form a category, where the composition of morphisms is given by the vertical composition. 

		Let $F_1: \calc \to \cald$, $F_2: \cald \to \mathcal{E}$ be $(G,\Gamma)$-crossed functors between $(G,\Gamma)$-crossed monoidal categories. By a standard argument, $F_2 F_1$ can be made into a $\Gamma$-graded monoidal functor.

		We make $F_2 F_1$ into a $(G,\Gamma)$-crossed monoidal functor by putting $\eta^{F_2 F_1}_g \coloneqq (\eta^{F_2}_g \ast \id_{F_1}) (\id_{F_2} \ast \eta_g^{F_1})$ for $g \in G$. Indeed, axiom 1 in Definition \ref{definition_crossed_functor} for $(F_2 F_1,\eta^{F_2 F_1})$ follows by a standard argument. Namely, for $g \in G$, let a crossing from $F \gamma^\calc(g)$ to $\gamma^\cald(g)F$ graphically denote $\eta_g^F$ in the 2-category of categories. Also, let a fork from $\gamma^\calc(g) \gamma^\cald(h)$ to $\gamma^\calc(g h)$ (resp. from $\gamma^\cald(g) \gamma^\cald(h)$ to $\gamma^\cald(g h)$) graphically denote $\chi^{\gamma^\calc}_{g,h}$ (resp. $\chi^{\gamma^\cald}_{g,h}$). Then, axiom 1 in Definition \ref{definition_crossed_functor} graphically means that a fork can pass under a crossing. Then, axiom 1 for $(F_2 F_1,\eta^{F_2 F_1})$ follows since a fork can pass under two crossings by a standard graphical calculation.

		Axiom 2 follows since
		\begin{align*}
			&(\eta^{F_2 F_1}_g)_{\lambda \mu} F_2 F_1 (J^{\gamma^\calc(g)}_{\lambda, \mu}) J^{F_2 F_1}_{{}^{\partial \mu \rhd_2 g} \lambda, {}^g \mu} \\
			&= (\eta^{F_2}_g)_{F_1(\lambda \mu)} F_2((\eta_g^{F_1})_{\lambda \mu} F_1(J^{\gamma^\calc(g)}_{\lambda, \mu})J^{F_1}_{{}^{\partial \mu \rhd_2 g} \lambda, {}^g \mu}) J^{F_2}_{F_1({}^{\partial \mu \rhd_2 g} \lambda), F_1({}^g \mu)} \\
			&= (\eta^{F_2}_g)_{F_1(\lambda \mu)} F_2({}^g J^{F_1}_{\lambda, \mu} J^{\gamma^\cald(g)}_{F_1(\lambda), F_1(\mu)}((\eta_{\partial \mu \rhd_2 g}^{F_1})_{\lambda} \otimes (\eta_g^{F_1})_{\mu})) J^{F_2}_{F_1({}^{\partial \mu \rhd_2 g} \lambda), F_1({}^g \mu)} \\
			&= {}^g F_2(J^{F_1}_{\lambda, \mu}) (\eta^{F_2}_g)_{F_1(\lambda)F_1(\mu)} F_2 ( J^{\gamma^\cald(g)}_{F_1(\lambda), F_1(\mu)}) J^{F_2}_{{}^{\partial \mu \rhd_2 g} F_1(\lambda), {}^g F_1(\mu)} (F_2((\eta_{\partial \mu \rhd_2 g}^{F_1})_{\lambda}) \otimes F_2((\eta_g^{F_1})_{\mu})) \\
			&= {}^g J^{F_2 F_1}_{\lambda, \mu} J^{\gamma^\mathcal{E}(g)}_{F_2 F_1(\lambda), F_2F_1(\mu)} ((\eta^{F_2 F_1}_{\partial \mu \rhd_2 g})_\lambda \otimes (\eta^{F_2 F_1}_{g})_\mu)
		\end{align*}
		for $\lambda \in \obj(\calc)$, $\mu \in \homog(\calc)$ and $g \in G$ by axiom 2 for $F_1$ and $F_2$ and the naturality of $\eta_g^{F_2}$ and $J^{F_2}$. Note that we used that $F_1$ is $\Gamma$-graded at the final equality. 
		
		It is easy to see that the horizontal composition of $(G,\Gamma)$-crossed monoidal transformation is again $(G,\Gamma)$-crossed monoidal, and therefore the composition of $(G,\Gamma)$-crossed monoidal functors is a bifunctor. 

		We have $\eta^{(F_3 F_2) F_1} = \eta^{F_3(F_2 F_1)}$ for composable $(G,\Gamma)$-crossed monoidal functors $F_1, F_2, F_3$ by definition. Moreover, we have $\eta^{\id_\cald F} = \eta^{F \id_\calc} = \eta^F$ for a $(G,\Gamma)$-crossed monoidal functor $F: \calc \to \cald$ by definition. Thus, we obtain a 2-category.  
	\end{proof}
\end{proposition}

\begin{definition}
	A \emph{$(G, \Gamma)$-crossed monoidal equivalence} between $(G,\Gamma)$-crossed monoidal categories is an equivalence in the 2-category of $(G,\Gamma)$-crossed monoidal categories (see e.g. \cite[Section 1.5]{MR2094071}). Equivalently, a $(G,\Gamma)$-crossed monoidal functor $F: \calc \to \cald$ between $(G,\Gamma)$-crossed monoidal categories is a \emph{$(G, \Gamma)$-crossed monoidal equivalence} if there exists a $(G, \Gamma)$-crossed monoidal functor $F^{-1}: \cald \to \calc$ with $(G,\Gamma)$-crossed natural isomorphisms $F^{-1}F \cong \id_\calc$ and $FF^{-1} \cong \id_\cald$.
\end{definition}

The following lemma is useful.

\begin{lemma}
	\label{lemma_equivalence}
	Let $F: \calc \to \cald$ be a $(G, \Gamma)$-crossed monoidal (resp. tensor) functor between $(G, \Gamma)$-crossed monoidal (resp. multiring) categories. Then, $F$ is a $(G, \Gamma)$-crossed monoidal equivalence if and only if it is an equivalence as a functor. 

	\begin{proof}
		We show the ``if'' part since the ``only if'' part is trivial. Suppose $F$ is an equivalence as a functor. Then, by a standard argument, an adjoint inverse $F^{-1}$ of $F$ is additive (resp. linear exact faithful).  Moreover, by a standard argument (see e.g. \cite[Remark 2.4.10]{egno}), we can make $F^{-1}$ into a monoidal functor and obtain a monoidal evaluation map $\mathrm{ev}^F: F^{-1}F \cong \id_\calc$ and a monoidal coevaluation map $\mathrm{coev}^F: \id_{\cald} \cong FF^{-1}$. 
		
		Since $\calc$ and $\cald$ are $\Gamma$-graded, for $\lambda \in \obj(\cald_s)$, we can consider the homogeneous decomposition $F^{-1}(\lambda) = \bigoplus_{t \in \Gamma} F^{-1}(\lambda)_t$ of $F^{-1}(\lambda)$. Since $F$ is $\Gamma$-graded and $FF^{-1} \cong \id_\cald$ by assumption, $FF^{-1}(\lambda)_t$ is a zero object unless $t = s$. Thus, $F^{-1}FF^{-1}(\lambda)_t$ is a zero object unless $t=s$, and therefore $F^{-1}$ is $\Gamma$-graded since $F^{-1} \cong F^{-1}FF^{-1}$.
		
		We make $F^{-1}$ into a $(G,\Gamma)$-crossed monoidal (resp. tensor) functor. Put 
		\begin{align*}
			\eta^{F^{-1}}_g \coloneqq (\mathrm{ev}^F \ast \id_{\gamma^{\calc}(g) F^{-1}}) (\id_{F^{-1}} \ast (\eta^F_{g})^{-1} \ast \id_{F^{-1}}) (\id_{F^{-1}\gamma^\cald(g)} \ast \mathrm{coev}^F) 
		\end{align*}
		for $g \in G$. We check the axioms in Definition \ref{definition_crossed_functor} for $(F^{-1}, \eta^{F^{-1}})$. Let us use the graphical notations in the proof of Proposition \ref{proposition_crossed_bicategory}. Then, $\eta^{F^{-1}}_g$ is graphically represented by a rotated crossing, and therefore axiom 1 for $(F^{-1}, \eta^{F^{-1}})$ follows from axiom 1 and a conjugate equation for $F$. 
	Axiom 2 follows since
		\begin{align*}
			&(\eta^{F^{-1}}_g)_{\lambda \mu} F^{-1}(J^{\gamma^\cald(g)}_{\lambda, \mu}) J^{F^{-1}}_{{}^{\partial \mu \rhd_2 g} \lambda, {}^g \mu}\\
			&= \mathrm{ev}^F_{{}^g F^{-1}(\lambda \mu)} F^{-1}((\eta^F_{g})^{-1}_{F^{-1} (\lambda \mu)} {}^g \mathrm{coev}^F_{\lambda \mu} J^{\gamma^\cald(g)}_{\lambda, \mu}) J^{F^{-1}}_{{}^{\partial \mu \rhd_2 g} \lambda, {}^g \mu} \\
			&= \mathrm{ev}^F_{{}^g F^{-1}(\lambda \mu)} \\
			&\quad F^{-1}((\eta^F_{g})^{-1}_{F^{-1} (\lambda \mu)}{}^g F(J^{F^{-1}}_{\lambda,\mu})  {}^g J^F_{F^{-1}(\lambda),F^{-1}(\mu)} J^{\gamma^\cald(g)}_{FF^{-1}(\lambda), FF^{-1}(\mu)} ({}^g \mathrm{coev}^F_{\lambda} \otimes {}^g \mathrm{coev}^F_{\mu}))J^{F^{-1}}_{{}^{\partial \mu \rhd_2 g} \lambda, {}^g \mu} \\
		&= {}^g J^{F^{-1}}_{\lambda,\mu} \mathrm{ev}^F_{{}^g (F^{-1}(\lambda)F^{-1}(\mu))}  \\
		&\quad F^{-1}((\eta^F_{g})^{-1}_{F^{-1} (\lambda) F^{-1} (\mu)} {}^g J^F_{F^{-1}(\lambda),F^{-1}(\mu)} J^{\gamma^\cald(g)}_{FF^{-1}(\lambda), FF^{-1}(\mu)} ({}^g \mathrm{coev}^F_{\lambda} \otimes {}^g \mathrm{coev}^F_{\mu}))J^{F^{-1}}_{{}^{\partial \mu \rhd_2 g} \lambda, {}^g \mu} \\
			&= {}^g J^{F^{-1}}_{\lambda,\mu} \mathrm{ev}^F_{{}^g (F^{-1}(\lambda)F^{-1}(\mu))} F^{-1}F(J^{\gamma^\calc(g)}_{F^{-1}(\lambda),F^{-1}(\mu)}) F^{-1}(J^F_{{}^{\partial \mu \rhd_2 g} \lambda, {}^g \mu}) \\
			&\quad F^{-1}((\eta^F_{\partial \mu \rhd_2 g})^{-1}_{F^{-1}(\lambda)} {}^g \mathrm{coev}^F_{\lambda} \otimes (\eta^F_{g})^{-1}_{F^{-1}(\mu)} {}^g \mathrm{coev}^F_{\mu}) J^{F^{-1}}_{{}^{\partial \mu \rhd_2 g} \lambda, {}^g \mu} \\ 
			&= {}^g J^{F^{-1}}_{\lambda,\mu} J^{\gamma^\calc(g)}_{F^{-1}(\lambda),F^{-1}(\mu)} \mathrm{ev}^F_{{}^g F^{-1}(\lambda){}^g F^{-1}(\mu)} J^{F^{-1}F}_{{}^{\partial \mu \rhd_2 g} \lambda, {}^g \mu} \\
			&\quad (F^{-1}((\eta^F_{\partial \mu \rhd_2 g})^{-1}_{F^{-1}(\lambda)} {}^g \mathrm{coev}^F_{\lambda}) \otimes F^{-1} ((\eta^F_{g})^{-1}_{F^{-1}(\mu)} {}^g \mathrm{coev}^F_{\mu}))  \\
			&= {}^g J^{F^{-1}}_{\lambda,\mu} J^{\gamma^\calc(g)}_{F^{-1}(\lambda),F^{-1}(\mu)} ((\eta^{F^{-1}}_{\partial \mu \rhd_2 g})_\lambda \otimes (\eta^{F^{-1}}_{g})_\mu)
		\end{align*}
		for $\lambda \in \obj(\cald)$, $\mu \in \homog(\cald)$ and $g \in G$ by the naturality of $J^{\gamma^\cald(g)}$, $J^{F^{-1}}$, $\mathrm{ev}^F$ and $\mathrm{coev}^F$, axiom 2 for $F$ and the monoidality of $\mathrm{ev}^F$ and $\mathrm{coev}^F$. Thus, $(F^{-1},\eta^{F^{-1}})$ is a $(G,\Gamma)$-crossed monoidal (resp. tensor) functor. 

		Finally, recall that $\eta^{F^{-1}F}_g = (\eta^{F^{-1}}_g \ast \id_F) (\id_{F^{-1}} \ast \eta_g^F)$ and $\eta^{FF^{-1}}_g = (\eta^{F}_g \ast \id_{F^{-1}}) (\id_{F} \ast \eta_g^{F^{-1}})$ for $g \in G$ by the proof of Proposition \ref{proposition_crossed_bicategory}. Then, $\mathrm{ev}^F$ and $\mathrm{coev}^F$ are $(G,\Gamma)$-crossed natural isomorphisms since
		\begin{align*}
			(\id_{\gamma^\cald(g)} \ast \mathrm{ev}^F) (\eta^{F^{-1}}_g \ast \id_F) (\id_{F^{-1}} \ast \eta_g^F) &= \mathrm{ev}^F \ast \id_{\gamma^\calc(g)} \\
			(\eta^{F}_g \ast \id_{F^{-1}}) (\id_{F} \ast \eta_g^{F^{-1}}) (\mathrm{coev}^F \ast \id_{\gamma^\calc(g)}) &= \id_{\gamma^\cald(g)} \ast \mathrm{coev}^F 
		\end{align*}
		for $g \in G$ by conjugate equations for $F$. Thus, $(F^{-1},\eta^{F^{-1}})$ is a $(G,\Gamma)$-crossed monoidal (resp. tensor) equivalence.
	\end{proof}
\end{lemma}

\begin{lemma}
	\label{lemma_transport_crossed}
	Let $\calc$ be a $(G,\Gamma)$-crossed monoidal category, and let $\cald$ be an additive monoidal category. Given an additive monoidal equivalence $F:\calc \to \cald$, we can endow $\cald$ a $(G,\Gamma)$-crossed structure so that $F$ is a $(G,\Gamma)$-crossed monoidal equivalence. 

	\begin{proof}
		First, we give a $\Gamma$-grading on $\cald$. For $s \in \Gamma$, put $\cald_s \coloneqq \calc_s$, where $\calc_s$ is regarded as an additive subcategory of $\cald$ by using $F$. Then, since $\calc$ is $\Gamma$-graded and $F$ is essentially surjective, any object of $\cald$ is a direct sum of homogeneous objects. Since $F$ is moreover an equivalence, $\obj(\cald_s) \cap \obj(\cald_t)$ consists of zero objects if $s \neq t$. Thus, $\cald = \bigoplus_{s \in \Gamma} \cald_s$. Note that $F$ is $\Gamma$-graded by the definition of the $\Gamma$-grading on $\cald$.
 
		Next, we give a $G$-action on $\cald$. Fix an adjoint inverse $F^{-1}$ and put $\gamma^\cald(g) \coloneqq F \gamma^\calc(g) F^{-1}$ for $g \in G$. By a standard argument, $\gamma^\cald$ defines a $G$-action on the category $\cald$ by putting $\chi^{\gamma^\cald}_{g,h} \coloneqq \id_{F \gamma^\calc(g)} \ast \mathrm{ev}^F \ast \id_{ \gamma^\calc(h) F^{-1}}$ and $\iota^{\gamma^\cald} \coloneqq \mathrm{coev}^F$. Moreover, axiom 1 in Definition \ref{definition_crossed_category} is satisfied since $F^{-1}$ is $\Gamma$-graded by the proof of Lemma \ref{lemma_equivalence}.

		Then, we give a $(G,\Gamma)$-crossed structure. Put $J^{\gamma^\cald(g)}_{\lambda, \mu} \coloneqq F(J^{\gamma^\calc(g)}_{F^{-1}(\lambda), F^{-1}(\mu)})$ and $\varphi^{\gamma^\cald(g)} \coloneqq F(\varphi^{\gamma^\calc(g)})$ for $\lambda \in \obj(\cald)$, $\mu \in \homog(\cald)$ and $g \in G$. Note that $J^{\gamma^\cald(g)}_{\lambda, \mu}$ is a well-defined isomorphism ${}^{\partial \mu \rhd_2 g} \lambda {}^g \mu \cong {}^g (\lambda \mu)$ since $F^{-1}$ is $\Gamma$-graded. Then, axioms 2 and 3 in Definition \ref{definition_crossed_category} follow from those for $\calc$. Thus, $\cald$ is a $(G,\Gamma)$-crossed monoidal category.

		Finally, we make $F$ into a $(G,\Gamma)$-crossed equivalence. By Lemma \ref{lemma_equivalence}, it suffices to make it into a $(G,\Gamma)$-crossed monoidal functor. We have already seen that $F$ is $\Gamma$-graded. Put $\eta^F_g \coloneqq \id_{F \gamma^\calc(g)} \ast (\mathrm{ev}^F)^{-1}$ for $g \in G$. It is easy to check that $(F,\eta^F)$ satisfies axiom 1 in Definition \ref{definition_crossed_functor}. Axiom 2 follow from the definition of $J^{\gamma^\cald}$ and $\varphi^{\gamma^\cald}$ and the monoidality of $\mathrm{ev}^F$.
	\end{proof}
\end{lemma}

\begin{definition}
	A $(G, \Gamma)$-crossed monoidal category $\calc$ is \emph{strictly $(G, \Gamma)$-crossed} if the action $\gamma^\calc$ is strict and $J^{\gamma^\calc}$ and $\varphi^{\gamma^\calc}$ in Definition \ref{definition_crossed_category} consist of identity maps.
\end{definition}

The following theorem should be called the coherence theorem for crossed tensor categories.

\begin{theorem}
	\label{theorem_strictification}
	Let $(G,\Gamma)$ be a matched pair of groups. Then, any $(G,\Gamma)$-crossed monoidal category is $(G,\Gamma)$-crossed equivalent to a strictly $(G,\Gamma)$-crossed strict monoidal category.

	\begin{proof}
		Recall from \cite[Subsection 2.2.3]{MR3076451} that for a monoidal category $\calc$, by letting $\st (\calc)$ be the strict monoidal category whose objects are strings of objects of $\calc$, we obtain a monoidal equivalence $\tilde{F}_\calc: \calc \simeq \st(\calc)$, which is defined by regarding the objects of $\calc$ as the strings of length 1. By definition, when $\calc$ is additive monoidal, $\st(\calc)$ is again additive monoidal and $\tilde{F}_\calc$ is additive. Therefore, $\st(\calc)$ can be made into a $(G,\Gamma)$-crossed monoidal category so that $\tilde{F}_\calc$ is a $(G,\Gamma)$-crossed monoidal equivalence by Lemma \ref{lemma_transport_crossed}. Thus, we may assume that $\calc$ is a strict monoidal category from the beginning.

		For a $(G,\Gamma)$-crossed strict monoidal category $\calc$, let us define the category $\calc(G)$ as in \cite[Section 4]{MR3671186}. Namely, an object of $\calc(G)$ is a pair $L=(L,\xi^L)$ of a family $L = \{ L_g \}_{g \in G}$ of objects of $\calc$ and a family of isomorphisms $\xi^L = \{ \xi_{g,h}^L: {}^g L_h \cong L_{g h} \}_{g,h \in G}$ such that $\xi_{gh,k}^L = \xi_{g,hk}^L {}^g \xi_{h,k}^L$ and $\xi_{e,h}^L = \id_{L_h}$ for all $g,h,k \in G$, where we suppressed some canonical morphisms (see Section \ref{section_preliminaries}). A morphism $f: L \to L'$ is a family $f = \{ f_g: L_g \to L_g' \}_{g \in G}$ of morphisms of $\calc$ such that $f_{gh} \xi^L_{g,h} = \xi^{L'}_{g,h} {}^g f_h$ for all $g,h \in G$.
		
		By a standard argument, $\calc(G)$ is an additive category. As in \cite[Section 4]{MR3671186}, we obtain a strict $G$-action $\gamma^{\calc(G)}$ on the additive category $\calc(G)$ by putting
		\begin{align*}
			\gamma^{\calc(G)}(g)(L) \coloneqq (\{ L_{h g} \}_h, \{ \xi^L_{k, hg} \}_{k,h})
		\end{align*}
		for $(L, \xi^L) \in \obj(\calc(G))$ and $g \in G$. Moreover, $\calc(G) = \bigoplus_{s \in \Gamma} \calc(G)_s$ by putting
		\begin{align*}
			\obj(\calc(G)_s) \coloneqq \{ L \in \obj(\calc(G)) \mid L_e \in \obj(\calc(G)_s) \}
		\end{align*}
		for $s \in \Gamma$ as in the proof of \cite[Corollary 5.3]{MR3671186}. For any $g \in G$ and $s\in \Gamma$, ${}^g \obj(\calc(G)_s) \subset \obj(\calc(G)_{g \rhd_1 s})$ since for $L \in \obj(\calc(G)_s)$ we have $({}^g L)_e = L_g \cong {}^g L_e \in \obj(\calc_{g \rhd_1 s})$. 

		We make $\calc(G)$ into a monoidal category by putting
		\begin{align*}
			L M \coloneqq (\{ L_{\partial M \rhd_2 h} M_h \}_h,\{ (\xi^L_{(h \rhd_1 \partial M) \rhd_2 g, \partial M \rhd_2 h} \otimes \xi^M_{g,h}) J^{\gamma^\calc(g)-1}_{L_{\partial M \rhd_2 h}, M_h} \}_{g,h})
		\end{align*}
		for $L \in \obj(\calc(G))$ and $M \in \homog(\calc(G))$. Note that $\xi^{LM}_{g,h}$ is indeed an isomorphism ${}^g (LM)_h \cong (LM)_{g h}$ by the matching equation $((h \rhd_1 \partial M) \rhd_2 g) (\partial M \rhd_2 h) = \partial M \rhd_2 (g h)$. By similar matching equations, we have 
		\begin{align*}
			\xi^L_{(k \rhd_1 \partial M) \rhd_2 (g h), \partial M \rhd_2 k} &= \xi^L_{((h k \rhd_1 \partial M) \rhd_2 g) ((k \rhd_1 \partial M) \rhd_2 h), \partial M \rhd_2 k} \\
			&= \xi^L_{(hk \rhd_1 \partial M) \rhd_2 g, ((k \rhd_1 \partial M) \rhd_2 h)(\partial M \rhd_2 k)} {}^{(hk \rhd_1 \partial M) \rhd_2 g} \xi^L_{(k \rhd_1 \partial M) \rhd_2 h, \partial M \rhd_2 k} \\
			&= \xi^L_{(hk \rhd_1 \partial M) \rhd_2 g, \partial M \rhd_2 h k} {}^{(hk \rhd_1 \partial M) \rhd_2 g} \xi^L_{(k \rhd_1 \partial M) \rhd_2 h, \partial M \rhd_2 k}
		\end{align*}
		and therefore
		\begin{align*}
			\xi^{LM}_{gh,k} &= (\xi^L_{(hk \rhd_1 \partial M) \rhd_2 g, \partial M \rhd_2 h k} {}^{(hk \rhd_1 \partial M) \rhd_2 g} \xi^L_{(k \rhd_1 \partial M) \rhd_2 h, \partial M \rhd_2 k} \otimes \xi^M_{g,hk} {}^g \xi^M_{h,k}) J^{\gamma^\calc(gh)-1}_{L_{\partial M \rhd_2 k}, M_k} \\
			&= (\xi^L_{(hk \rhd_1 \partial M) \rhd_2 g, \partial M \rhd_2 h k} {}^{(hk \rhd_1 \partial M) \rhd_2 g} \xi^L_{(k \rhd_1 \partial M) \rhd_2 h, \partial M \rhd_2 k} \otimes \xi^M_{g,hk} {}^g \xi^M_{h,k}) \\
			&\qquad J^{\gamma^\calc(g)-1}_{{}^{(k \rhd_1 \partial M) \rhd_2 h} L_{\partial M \rhd_2 k}, {}^h M_k} {}^g J^{\gamma^\calc(h)-1}_{L_{\partial M \rhd_2 k}, M_k} \\
			&= (\xi^L_{(hk \rhd_1 \partial M) \rhd_2 g, \partial M \rhd_2 h k} \otimes \xi^M_{g,hk}) J^{\gamma^\calc(g)-1}_{L_{\partial M \rhd_2 h k} , M_{h k}} {}^g ((\xi^L_{(k \rhd_1 \partial M) \rhd_2 h, \partial M \rhd_2 k} \otimes \xi^M_{h,k}) J^{\gamma^\calc(h)-1}_{L_{\partial M \rhd_2 k}, M_k}) \\
			&= \xi^{LM}_{g,hk} {}^g \xi^{LM}_{h,k}
		\end{align*}
		for $g,h,k \in G$ by axiom 3 in Definition \ref{definition_crossed_category} and the naturality of $J^{\gamma^\calc(g)}$. We also have $\xi^{LM}_{e,h} = (\xi^L_{e, \partial M \rhd_2 h} \otimes \xi^M_{e,h}) = \id$ for $g \in G$ by axiom 3 in Definition \ref{definition_crossed_category} and $(h \rhd_1 \partial M) \rhd_2 e = e$. Thus, $LM$ is well-defined. Moreover, for morphisms $f:L \to L'$ and $\tilde{f}: M \to M'$ of $\calc(G)$, define a morphism $f \otimes \tilde{f} : LM \to L' M'$ by putting $(f \otimes \tilde{f})_g \coloneqq f_{\partial M \rhd_2 g} \otimes \tilde{f}_g$ for $g \in G$. Then,
		\begin{align*}
			(f \otimes \tilde{f})_{gh} \xi^{LM}_{g,h} &= (f_{\partial M \rhd_2 gh} \xi^L_{(h \rhd_1 \partial M) \rhd_2 g, \partial M \rhd_2 h} \otimes \tilde{f}_{gh} \xi^M_{g,h} ) J^{\gamma^\calc(g)-1}_{L_{\partial M \rhd_2 h,M_h}} \\
			&= (\xi^{L'}_{(h \rhd_1 \partial M') \rhd_2 g, \partial M' \rhd_2 h} {}^{(h \rhd_1 \partial M) \rhd_2 g} f_h \otimes \xi^{M'}_{g,h} {}^g \tilde{f}_h) J^{\gamma^\calc(g)-1}_{L_{\partial M \rhd_2 h,M_h}} \\
			&= (\xi^{L'}_{(h \rhd_1 \partial M') \rhd_2 g, \partial M' \rhd_2 h} \otimes \xi^{M'}_{g,h}) J^{\gamma^\calc(g)-1}_{L'_{\partial M' \rhd_2 h,M'_h}} {}^g (f \otimes \tilde{f})
		\end{align*}
		for $g,h \in G$ by the naturality of $J^{\gamma^\calc(g)}$. Note that we may assume $\partial M = \partial M'$ since otherwise $\tilde{f} = 0$. Therefore, $f \otimes \tilde{f}$ is well-defined. Thus, the product $\otimes$ defines a biadditive bifunctor $\calc(G) \times \calc(G) \to \calc(G)$. Moreover, it makes $\calc(G)$ into a $\Gamma$-graded strict monoidal category with unit $\mathbf{1}_{\calc(G)} \coloneqq (\{ \mathbf{1}_{\calc} \}_{h}, \{ \varphi^{\gamma^\calc(g)-1} \}_{g,h})$. Indeed,
		\begin{align*}
			((LM)N)_h = L_{\partial M \rhd_2 (\partial N \rhd_2 h)} M_{\partial N \rhd_2 h} N_h = L_{\partial (MN) \rhd_2 h} M_{\partial N \rhd_2 h} N_h = (L(MN))_h
		\end{align*}
		for $L \in \obj(\calc(G))$, $M,N \in \homog(\calc(G))$ and $h \in G$ since $\rhd_2$ is an action and $\partial M \partial N = \partial M_e \partial N_e = \partial (M_e N_e) = \partial (M_{\partial N \rhd_2 e} N_e) = \partial (M N)$. Moreover, since
		\begin{align*}
			\xi^{(LM)N}_{g,h} &= (\xi^{LM}_{(h \rhd_1 \partial N) \rhd_2 g, \partial N \rhd_2 h} \otimes \xi^N_{g,h}) J^{\gamma^\calc(g)-1}_{(LM)_{\partial N \rhd_2 h}, N_h} \\
			&= (\xi^L_{((\partial N \rhd_2 h)\rhd_1 \partial M) \rhd_2 ((h \rhd_1 \partial N) \rhd_2 g),\partial M \partial N \rhd_2 h} \otimes \xi^{M}_{(h \rhd_1 \partial N) \rhd_2 g, \partial N \rhd_2 h} \otimes \xi^N_{g,h}) \\
			&\qquad (J^{\gamma^\calc(g)-1}_{L_{\partial M \partial N \rhd_2 h}, M_{\partial N \rhd_2 h}} \otimes \id_{{}^g N_h}) J^{\gamma^\calc(g)-1}_{L_{\partial M \partial N \rhd_2 h} M_{\partial N \rhd_2 h}, N_h} \\
			&= (\xi^L_{(h \rhd_1 \partial M \partial N)\rhd_2 g,\partial M \partial N \rhd_2 h} \otimes \xi^{M}_{(h \rhd_1 \partial N) \rhd_2 g, \partial N \rhd_2 h} \otimes \xi^N_{g,h}) \\
			&\qquad (J^{\gamma^\calc(g)-1}_{L_{\partial M \partial N \rhd_2 h}, M_{\partial N \rhd_2 h}} \otimes \id_{{}^g N_h}) J^{\gamma^\calc(g)-1}_{L_{\partial M \partial N \rhd_2 h} M_{\partial N \rhd_2 h}, N_h} \\
			\xi^{L(MN)}_{g,h} &= (\xi^L_{(h \rhd_1 \partial M \partial N)\rhd_2 g,\partial M \partial N \rhd_2 h} \otimes \xi^{MN}_{g,h}) J^{\gamma^\calc(g)-1}_{L_{\partial M \partial N \rhd_2 h}, (MN)_h} \\
			&= (\xi^L_{(h \rhd_1 \partial M \partial N)\rhd_2 g,\partial M \partial N \rhd_2 h} \otimes \xi^{M}_{(h \rhd_1 \partial N) \rhd_2 g, \partial N \rhd_2 h} \otimes \xi^N_{g,h}) \\
			&\qquad (\id_{{}^{(h \rhd_1 \partial M \partial N) \rhd_2 g} L_{\partial M \partial N \rhd_2 h}} \otimes J^{\gamma^\calc(g)-1}_{M_{\partial N \rhd_2 h}, N_h})J^{\gamma^\calc(g)-1}_{L_{\partial M \partial N \rhd_2 h}, M_{\partial N \rhd_2 h} N_h}
		\end{align*}
		for $g,h \in G$ by the matching relation $h \rhd_1 \partial M \partial N = ((\partial N \rhd_2 h) \rhd_1 \partial M)(h \rhd_1 \partial N)$, we obtain $\xi^{(LM)N} = \xi^{L(MN)}$ by axiom 2 in Definition \ref{definition_crossed_category}. The unit $\mathbf{1}_{\calc(G)}$ is indeed an object of $\calc(G)$ by axiom 3 in Definition \ref{definition_crossed_category}. We have $(\mathbf{1}_{\calc(G)} M)_h = M_h = (M \mathbf{1}_{\calc(G)})_h$ since $\partial \mathbf{1}_\calc = e$. Moreover, $\xi^{\mathbf{1}_{\calc(G)}M} = \xi^{M \mathbf{1}_{\calc(G)}} = \xi^{M}$ follows from axiom 2 in Definition \ref{definition_crossed_category}.
		
		For $L \in \obj(\calc(G))$, $M \in \homog(\calc(G))$ and $g,h \in G$, the matching relation $\partial M \rhd_2 hg = ((g \rhd_1 \partial M) \rhd_2 h) (\partial M \rhd_2 g)$, we have
		\begin{align*}
			({}^g (LM))_h &= (LM)_{hg} = L_{\partial M \rhd_2 hg} M_{hg} = L_{((g \rhd_1 \partial M)\rhd_2 h)(\partial M \rhd_2 g)} M_{hg} \\
			&= ({}^{\partial M \rhd_2 g} L)_{(g \rhd_1 \partial M)\rhd_2 h} ({}^g M)_h = ({}^{\partial M \rhd_2 g} L)_{\partial {}^g M \rhd_2 h} ({}^g M)_h \\
			&= ({}^{\partial M \rhd_2 g} L {}^g M)_h.
		\end{align*}
		Moreover,
		\begin{align*}
			\xi^{{}^{\partial M \rhd_2 g} L {}^g M}_{k,h} &= (\xi^{{}^{\partial M \rhd_2 g} L}_{(h \rhd_1 \partial {}^g M) \rhd_2 k, \partial {}^g M \rhd_2 h} \otimes \xi^{{}^g M}_{k,h}) J^{\gamma^\calc(k)-1}_{({}^{\partial M \rhd_2 g} L)_{\partial {}^g M \rhd_2 h},({}^g M)_h} \\
			&= (\xi^{{}^{\partial M \rhd_2 g} L}_{(hg \rhd_1 \partial M) \rhd_2 k, (g \rhd_1 \partial M) \rhd_2 h} \otimes \xi^{{}^g M}_{k,h}) J^{\gamma^\calc(k)-1}_{({}^{\partial M \rhd_2 g} L)_{(g \rhd_1 \partial M) \rhd_2 h},({}^g M)_h} \\
			&= (\xi^{L}_{(hg \rhd_1 \partial M) \rhd_2 k, \partial M \rhd_2 hg} \otimes \xi^{M}_{k,hg}) J^{\gamma^\calc(k)-1}_{L_{\partial M \rhd_2 hg},M_{hg}} \\
			&= \xi^{LM}_{k,hg} = \xi^{{}^g(LM)}_{k,h}
		\end{align*}
		for $g,h,k \in G$. Therefore, ${}^g(LM) = {}^{\partial M \rhd_2 g} L {}^g M$. We also have ${}^g \mathbf{1}_{\calc(G)} = \mathbf{1}_{\calc(G)}$ by definition. Thus, the $G$-action $\gamma^{\calc(G)}$ and the $\Gamma$-grading $\calc(G) = \bigoplus_{s \in \Gamma} \calc(G)_s$ make $\calc(G)$ into a strictly $(G,\Gamma)$-crossed strict monoidal category. 

		The functor $I_\calc: \calc \to \calc(G); \lambda \mapsto (\{ {}^h \lambda \}_h, \{ (\chi^{\gamma^\calc}_{g,h})_\lambda \}_{g,h})$ is an equivalence of categories as in the proof of \cite[Proposition 4.1]{MR3671186}. It is also $\Gamma$-graded by definition. We can make it into a monoidal functor by putting $J^{I_\calc}_{\lambda, \mu} \coloneqq \{ J^{\gamma^\calc(h)}_{\lambda,\mu} \}_h$ for $\lambda \in \obj(\calc)$ and $\mu \in \homog(\calc)$ and $\varphi^{I_\calc} \coloneqq \{ \varphi^{\gamma^\calc(h)} \}_h$. Indeed, they are morphisms of $\calc(G)$ by axiom 3 in Definition \ref{definition_crossed_category} since $\xi^{I_\calc(\lambda) I_\calc(\mu)}_{g,h} = J^{\gamma^\calc(g)-1}_{{}^{\partial \mu \rhd_2 h} \lambda,{}^h \mu}$ and $\xi^{\mathbf{1}_{\calc(G)}}_{g,h} = \varphi^{\gamma^\calc(g)-1}$ by definition. They satisfy the coherence relations for a monoidal functor by axiom 2 in Definition \ref{definition_crossed_category}. 

		Moreover, we can make $I_\calc$ into a $(G,\Gamma)$-crossed monoidal functor and therefore prove the theorem by Lemma \ref{lemma_equivalence}. Indeed, since
		\begin{align*}
			I_\calc ({}^g \lambda) &= (\{ {}^{h} ({}^g \lambda) \}_h,\{ (\chi^{\gamma^\calc}_{k,h})_{{}^g \lambda} \}_{k,h}) \\
			{}^g I_\calc(\lambda) &= (\{ {}^{h g} \lambda \}_h,\{ (\chi^{\gamma^\calc}_{k,hg})_\lambda \}_{k,h})
		\end{align*}
		for $\lambda \in \obj(\calc)$ and $g \in G$ by definition, we can put $(\eta^{I_\calc}_g)_\lambda \coloneqq \{(\chi^{\gamma^\calc}_{h,g})_\lambda\}_h : I_\calc ({}^g \lambda) \cong {}^g I_\calc(\lambda)$, which is indeed a morphism of $\calc(G)$ by coherence. Then, the first two axioms in Definition \ref{definition_crossed_functor} are satisfied by coherence. The remaining two axioms follow by coherence from the definition of $J^{I_\calc}$ and $\varphi^{I_\calc}$, respectively, and axiom 3 in Definition \ref{definition_crossed_category}. Thus, $(I_\calc, \eta^{I_\calc})$ is a $(G,\Gamma)$-crossed monoidal functor.
		\end{proof}
\end{theorem}

\begin{corollary}
	\label{corollary_coherence}
	Let $\calc$ be a $(G,\Gamma)$-crossed monoidal category. Define a set $W = \sqcup_{n \ge 0} W_n$ of words recursively by the following rules: $\mathbf{1} \in W_0$, $- \in W_1$, $\otimes \in W_2$, $\gamma(g) \in W_1$ for $g \in G$ and $w'((w_i)_i) \in W_{\sum_i n_i}$ for $w' \in W_n$ and a tuple $(w_i)_{i=1}^n$ with $w_i \in W_{n_i}$. For $w \in W_n$, let $w$ again denote the functor $\calc^n \to \calc$ corresponding to $w$, where $\calc^0$ denotes the category with only one object and its identity morphism. Define a set $I$ of morphisms recursively by the following rules: the components of $a^\calc,l^\calc,r^\calc, J^{\gamma^\calc(g)}$, $\chi^{\gamma^\calc}_{g,h}$ and $\iota^{\gamma^\calc}$ are in $I$ for any $g,h \in G$, $\varphi^{\gamma^\calc(g)} \in I$ for any $g \in G$, $f^{-1} \in I$ for any $f \in I$, $f \circ f' \in I$ for any composable $f, f' \in I$ and $w((f_i)_i)$ for any $(f_i)_{i=1}^n \subset I$ and $w \in W_n$. Then, for any $w, w' \in W_n$ and objects $\{\lambda_i\}_{i=1}^n$ of $\calc$, a morphism $w((\lambda_i)_i) \to w'((\lambda_i)_i)$ in $I$ is unique if it exists. 

	\begin{proof}
		By theorem \ref{theorem_strictification}, there exists a $(G,\Gamma)$-crossed monoidal equivalence $F$ from $\calc$ to a strictly $(G,\Gamma)$-crossed strict monoidal category $\tilde{\calc}$. For $w \in W_n$ and $\{ \lambda_i \}_{i=1}^n \subset \obj(\calc)$, we define a natural isomorphism $C^w: F w \cong \tilde{w} F^n$, where $\tilde{w}$ denotes the functor $\tilde{\calc}^n \to \tilde{\calc}$ corresponding to $w$, recursively by the following rules: $C^{\mathbf{1}} \coloneqq \varphi^{F-1}$, $C^{\otimes} \coloneqq J^{F-1}$, $C^{\gamma(g)} \coloneqq \eta^F_g$ for $g \in G$ and $C^{v((w_i)_i)} \coloneqq (\id_{\tilde{v}} \ast (C^{w_i})_i) (C^v \ast \id_{(w_i)_i})$ for $\{ w_i \}_{i=1}^n \subset W$ and $v \in W_n$. 

		We show that $C^w_{(\lambda_i)_i} = C^{w'}_{(\lambda_i)_i} F(f)$ for any $f: w((\lambda_i)_i) \to w'((\lambda_i)_i)$ in $I$, which uniquely determines $f$. It suffices to prove this when $f$ is a generator of $I$. The statement for $a^\calc, l^\calc, r^\calc$ follows since $F$ is monoidal, which is standard. The statement for $\chi^{\gamma^\calc}, \iota^{\gamma^\calc}, J^{\gamma^\calc(g)}$ and $\varphi^{\gamma^\calc(g)}$ follows from the axioms for $F = (F,\eta^F)$ in Definition \ref{definition_crossed_functor}.
	\end{proof}
\end{corollary}

\section{Crossed center}
\label{section_center}

In this section, we give the group-crossed version of the Drinfeld center construction (Theorem \ref{theorem_crossed_center}) by generalizing the graded center construction (see Example \ref{example_graded_center}) by Turaev and Virelizier \cite{MR3053524}. 

\begin{definition}
	Let $\calc$ be a $(G,\Gamma)$-crossed monoidal category. Then, the \emph{$(G,\Gamma)$-crossed center} $\rmz^G_\Gamma(\calc)$ of $\calc$ is defined as follows. For $g \in G$, an object of $\rmz^G_\Gamma(\calc)_g$ is the pair $(\lambda, h^\lambda)$ of $\lambda \in \obj(\calc)$ and a natural isomorphism $h^\lambda = \{ h^\lambda_\nu: \lambda \nu \cong {}^g \nu \lambda \}_{\nu \in \obj(\calc_e)}$ such that $h^\lambda_{\nu_1 \nu_2} = (\id_{{}^g \nu_1} \otimes h^\lambda_{\nu_2}) (h^\lambda_{\nu_1} \otimes \id_{\nu_2})$ for all $\nu_1, \nu_2 \in \obj(\calc_e)$. A morphism $f: (\lambda, h^\lambda) \to (\mu, h^\mu)$ of $\rmz^G_\Gamma(\calc)_g$ is a morphism $f: \lambda \to \mu$ of $\calc$ such that $h^\mu_\nu (f \otimes \id_\nu) = (\id_{{}^g \nu} \otimes f) h^\lambda_\nu$ for all $\nu \in \obj(\calc_e)$. Then, $\rmz^G_\Gamma(\calc)_g$ is an additive category, where the composition of morphisms and identity morphisms are given by those of $\calc$. The category $\rmz^G_\Gamma(\calc)$ is defined to be $\bigoplus_{g \in G} \rmz^G_\Gamma(\calc)_g$.
\end{definition}

\begin{lemma}
	Let $\calc$ be a $(G,\Gamma)$-crossed monoidal category. Then, the $(G,\Gamma)$-crossed center $\rmz^G_\Gamma(\calc)$ of $\calc$ becomes a $G \times \Gamma$-graded monoidal category by putting 
	\begin{align*}
		\obj(\rmz^G_\Gamma(\calc)_{(g,s)}) \coloneqq \{ (\lambda, h^\lambda) \in \obj(\rmz^G_\Gamma(\calc)_g) \mid \lambda \in \obj(\calc_s) \}
	\end{align*}
	for $(g,s) \in G \times \Gamma$, $(\lambda, h^\lambda) \otimes (\mu, h^\mu) \coloneqq (\lambda \mu, (h^\lambda \otimes \id)(\id \otimes h^\mu))$ for $(\lambda, h^\lambda), (\mu, h^\mu) \in \homog(\rmz^G_\Gamma(\calc))$ and $\mathbf{1}_{\rmz^G_\Gamma(\calc)} \coloneqq (\mathbf{1}_\calc, \id)$. 

	\begin{proof}
		By coherence (Corollary \ref{corollary_coherence}), we can do the same graphical calculation as that for the ordinary center of $\calc$.
	\end{proof}
\end{lemma}

\begin{lemma}
	\label{lemma_crossed_dual}
	Let $\calc$ be a $(G,\Gamma)$-crossed monoidal category. If $\lambda \in \homog(\calc)$ has a left dual $\lambda^\vee$, then, for any $g \in G$, ${}^g \lambda$ has a left dual ${}^{\partial \lambda \rhd_2 g} (\lambda^\vee)$ with evaluation and coevaluation maps given by $\mathrm{ev}_{{}^g \lambda} \coloneqq {}^g \mathrm{ev}_\lambda$ and $\mathrm{coev}_{{}^g \lambda} \coloneqq {}^{\partial \lambda \rhd_2 g} \mathrm{coev}_{\lambda}$, where we suppressed canonical isomorphisms in Corollary \ref{corollary_coherence}. We also have a similar statement for right duals.

	\begin{proof}
		We only show the statement for left duals. Note that the target of $\mathrm{coev}_{{}^g \lambda}$ is ${}^{\partial \lambda \rhd_2 g} (\lambda \lambda^\vee) = {}^{\partial \lambda^\vee \rhd_2 (\partial \lambda \rhd_2 g)} \lambda {}^{\partial \lambda \rhd_2 g} (\lambda^\vee) = {}^g \lambda {}^{\partial \lambda \rhd_2 g} (\lambda^\vee)$. Then, we have
		\begin{align*}
			(\id_{{}^g \lambda} \otimes {}^g \mathrm{ev}_\lambda) ({}^{\partial \lambda \rhd_2 g} \mathrm{coev}_{\lambda} \otimes \id_{{}^g \lambda}) &=  ({}^{e \rhd_2 g} \id_{\lambda} \otimes {}^g \mathrm{ev}_\lambda) ({}^{\partial \lambda \rhd_2 g} \mathrm{coev}_{\lambda} \otimes {}^g \id_{\lambda}) \\
			&= {}^g ((\id_{\lambda} \otimes \mathrm{ev}_\lambda) (\mathrm{coev}_{\lambda} \otimes \id_{\lambda})) = \id_{{}^g \lambda}
		\end{align*}
		and
		\begin{align*}
			&({}^g \mathrm{ev}_\lambda \otimes \id_{{}^{\partial \lambda \rhd_2 g}(\lambda^\vee)}) (\id_{{}^{\partial \lambda \rhd_2 g}(\lambda^\vee)} {}^{\partial \lambda \rhd_2 g} \mathrm{coev}_{\lambda}) \\
			&= ({}^{\partial \lambda^\vee \rhd_2(\partial \lambda \rhd_2 g)} \mathrm{ev}_\lambda \otimes {}^{\partial \lambda \rhd_2 g} \id_{\lambda^\vee}) ({}^{e \rhd_2(\partial \lambda \rhd_2 g)} \id_{\lambda^\vee} {}^{\partial \lambda \rhd_2 g} \mathrm{coev}_{\lambda}) \\
			&= {}^{\partial \lambda \rhd_2 g}((\mathrm{ev}_\lambda \otimes \id_{\lambda^\vee}) (\id_{\lambda^\vee} \otimes \mathrm{coev}_{\lambda})) = \id_{{}^{\partial \lambda \rhd_2 g} \lambda^\vee}
		\end{align*}
		by coherence (Corollary \ref{corollary_coherence}).
	\end{proof}
\end{lemma}

\begin{lemma}
	Let $\calc$ be a $(G,\Gamma)$-crossed monoidal category. If $\calc$ is rigid, then so is $\rmz^G_\Gamma(\calc)$. 

	\begin{proof}
		We can rotate a crossing and do the usual graphical calculations by Lemma \ref{lemma_crossed_dual} and coherence (Corollary \ref{corollary_coherence}). 
	\end{proof}
\end{lemma}

Then, we give crossed structures.

\begin{lemma}
	\label{lemma_center_g_action}
	The $(G,\Gamma)$-crossed center $\rmz^G_\Gamma(\calc)$ of a $(G,\Gamma)$-crossed monoidal category $\calc$ becomes a $(G,G \times \Gamma)$-crossed monoidal category, where $(G, G \times \Gamma)$ denotes the matched pair $(G, G \times \Gamma, \rhd_1^G, \rhd_2^G)$ defined by $g \rhd_1^G (h,t) \coloneqq ((t \rhd_2 g)hg^{-1},g \rhd_1 t)$ and $(h,t) \rhd_2 g \coloneqq t \rhd_2 g$, by putting ${}^g (\lambda, h^\lambda) \coloneqq ({}^g \lambda, h^{{}^g \lambda})$ and $h^{{}^g \lambda}_{\nu} \coloneqq {}^g h^\lambda_{{}^{g^{-1}} \nu} $ for $\nu \in \obj(\calc_e)$.

	\begin{proof}
		Immediate from coherence (Corollary \ref{corollary_coherence}).
	\end{proof}
\end{lemma}

We define a \emph{$(G,\Gamma)$-crossed pivotal tensor category} to be a $(G,\Gamma)$-crossed tensor category $\calc$ with ${}^g \delta^\calc_{\lambda} = \delta^\calc_{{}^g \lambda}$ for any $\lambda \in \obj(\calc)$. Note that this is meaningful thanks to Lemma \ref{lemma_crossed_dual}.

The following lemma is the group-crossed version of \cite[Section 4.4]{MR3053524}. Recall that a pivotal $\Gamma$-graded tensor category $\calc$ is \emph{non-singular} if for any $s \in \Gamma$, $\calc_s$ has at least one object with nonzero dimension. 

\begin{lemma}
	\label{lemma_center_gamma_action}
	The $(G,\Gamma)$-crossed center $\rmz^{(G,\Gamma)}(\calc)$ of a non-singular $(G,\Gamma)$-crossed pivotal tensor category $\calc$ admits a $(\Gamma,G \times \Gamma)$-crossed structure, where $(\Gamma,G \times \Gamma)$ denotes the matched pair $(\Gamma,G \times \Gamma, \rhd_1^\Gamma ,\rhd_2^\Gamma)$ defined by $s \rhd_1^\Gamma (h,t) \coloneqq (s \rhd_2 h, (h \rhd_1 s)ts^{-1})$ and $(h,t) \rhd_2^\Gamma s \coloneqq h \rhd_1 s$.

	\begin{proof}
		For $\lambda = (\lambda,h^\lambda) \in \rmz^{(G,\Gamma)}(\calc)_{(h,t)}$ and $\mu \in \obj(\calc_s)$ with dimension $d_\mu \neq 0$, put
		\begin{align*}
			e_{\lambda}^\mu \coloneqq d_\mu^{-1} (\mathrm{ev}_{{}^{s \rhd_2 h} (\mu^\vee)} \otimes \id_{{}^h \mu \lambda \mu^\vee}) (\delta^\calc_{{}^h \mu} \otimes h^\lambda_{\mu^\vee \mu} \otimes \id_{\mu^\vee}) (\id_{{}^h \mu \lambda \mu^\vee} \otimes \mathrm{coev}_{\mu}),
		\end{align*}
		which is well-defined since $({}^h \mu)^\vee = {}^{s \rhd_2 h} (\mu^\vee)$ by Lemma \ref{lemma_crossed_dual}. 

		Then, $e^\mu_{\lambda}$ is an idempotent in $\en_\calc({}^h \mu \lambda \mu^\vee)$ by the same argument as in the proof of \cite[Lemma 4.2]{MR3053524}. Let $(E_{\lambda}^\mu,s^\mu_\lambda,r^\mu_\lambda)$ denote a retract of $e^\mu_{\lambda}$ and put
		\begin{align*}
			h^{E^\mu_\lambda}_\nu &\coloneqq d_\mu^{-1} (\id_{{}^{s \rhd_2 h} \nu} \otimes r^\mu_\lambda) (\mathrm{ev}_{{}^{s \rhd_2 h} (\mu^\vee)} \otimes \id_{{}^{s \rhd_2 h} \nu {}^h \mu \lambda \mu^\vee}) (\delta^\calc_{{}^h \mu} \otimes h^\lambda_{\mu^\vee \nu \mu} \otimes \id_{\mu^\vee}) (\id_{{}^h \mu \lambda \mu^\vee \nu} \otimes \mathrm{coev}_{\mu}) \\
			&\qquad (s^\mu_\lambda \otimes \id_\nu)
		\end{align*}
		for $\nu \in \obj(\calc_e)$. Then, $h^{E^\mu_\lambda}_{\mathbf{1}_\calc} = \id_{E^\mu_\lambda}$ by definition. Moreover, by the same argument as in the proof of \cite[Lemma 4.3]{MR3053524}, we obtain $h^{E^\mu_\lambda}_{\nu_1 \nu_2} = (\id_{{}^h \nu_1} \otimes h^{E^\mu_\lambda}_{\nu_2})(h^{E^\mu_\lambda}_{\nu_1} \otimes \id_{\nu_2})$ and therefore
		\begin{align*}
			\gamma^\Gamma_\mu(\lambda) \coloneqq (E^\mu_\lambda,h^{E^\mu_\lambda}) \in \obj(\rmz^{(G,\Gamma)}(\calc)_{(s \rhd_2 h,(h \rhd_1 s) t s^{-1})}).
		\end{align*}

		By putting $\gamma^\Gamma_\mu(f) \coloneqq r_{\lambda'}^\mu f s_\lambda^\mu$ for $f \in \hom_{\rmz^{(G,\Gamma)}(\calc)_{(h,t)}}(\lambda,\lambda')$, we can regard $\gamma^\Gamma_\mu$ as an endofunctor of $\calc$ as in the proof of \cite[Lemma 4.3]{MR3053524}. Then, the coherence theorem (Corollary \ref{corollary_coherence}) makes graphical calculation possible, and the same argument as in the proof of \cite[Lemmata 4.3, 4.4 and 4.5]{MR3053524} shows that for a family $\{ \zeta_s \}_{s \in \Gamma}$ of homogeneous objects with nonzero dimensions, the assignment $s \mapsto \gamma^\Gamma(s) \coloneqq \gamma^\Gamma_{\zeta_s}$ defines a $(\Gamma, G \times \Gamma)$-crossed structure. 
	\end{proof}
\end{lemma}

By a tedious calculation, it can be checked that $(G \bowtie \Gamma, G \times \Gamma, \tilde{\rhd}_1, \tilde{\rhd}_2)$, where $(g,s) \tilde{\rhd}_1 (h,t) \coloneqq g \rhd_1^G (s \rhd_1^\Gamma (h,t))$ and $(h,t) \tilde{\rhd}_2 (g,s) \coloneqq ((s \rhd_1^\Gamma h) \rhd_2^G g, h \rhd_2^\Gamma s)$, is a matched pair. It is easy to check that the homomorphisms $\phi, \psi: G \times \Gamma \to G \bowtie \Gamma$ defined by putting $\phi(h,t) \coloneqq (e,t)$ and $\psi(h,t) \coloneqq (h,t)$ for $(h,t) \in G \times \Gamma$ give a braiding on the matched pair $(G \bowtie \Gamma, G \times \Gamma)$.

We combine two crossed structures. The following lemma is the key. 

\begin{lemma}
	\label{lemma_bicross_action}
	Let $\calc$ be a $G \times \Gamma$-graded monoidal category with $(G, G \times \Gamma)$-crossed and $(\Gamma, G \times \Gamma)$-crossed structures. Let $\gamma^G$ and $\gamma^\Gamma$ denote, respectively, the actions of $G$ and $\Gamma$ on $\calc$ given by the crossed structures. Suppose we are given a family of natural isomorphisms
	\begin{align*}
		\{ \sigma_{g,s} : \gamma^\Gamma(s) \gamma^G(g) \cong \gamma^G((s \rhd_2 g^{-1})^{-1}) \gamma^\Gamma(g^{-1} \rhd_1 s) \}_{g \in G, s \in \Gamma}
	\end{align*}
	such that
	\begin{align*}
		&(\sigma_{g,s})_{\lambda \mu} ({}^s J^{\gamma^G(g)}_{\lambda, \mu}) J^{\gamma^\Gamma(s)}_{{}^{\partial \mu \rhd_2^G g} \lambda, {}^g \mu} \\
		&= {}^{(s \rhd_2 g^{-1})^{-1}}J^{\gamma^\Gamma(g^{-1} \rhd_1 s)}_{\lambda, \mu}  J^{\gamma^G((s \rhd_2 g^{-1})^{-1})}_{{}^{\partial \mu \rhd_2^\Gamma (g^{-1} \rhd_1 s)} \lambda, {}^{g^{-1} \rhd_1 s} \mu} ((\sigma_{\partial \mu \rhd_2^G g,(g \rhd_1^G \partial \mu) \rhd_2^\Gamma s})_\lambda \otimes (\sigma_{g,s})_\mu) \\
		&(\sigma_{g,s})_{\mathbf{1}_\calc} ({}^s \varphi^{\gamma^G(g)}) \varphi^{\gamma^\Gamma(s)} = {}^{(s \rhd_2 g^{-1})^{-1}} \varphi^{\gamma^\Gamma(g^{-1} \rhd_1 s)} \varphi^{\gamma^G((s \rhd_2 g^{-1})^{-1})}
	\end{align*}
	and
	\begin{align*}
		\sigma_{g, ss'}(\chi^{\gamma^\Gamma}_{s, s'} \ast \id_{\gamma^G(g)}) &= (\id \ast \chi^{\gamma^\Gamma}_{(s' \rhd_2 g^{-1})^{-1} \rhd_1 s, g \rhd_1 s'}) (\sigma_{(s' \rhd_2 g^{-1})^{-1}, s} \ast \id_{\gamma^\Gamma(g^{-1} \rhd_1 s')}) (\id_{\gamma^\Gamma(s)} \ast \sigma^{g, s'}) \\
		\sigma_{gg', s} (\id_{\gamma^\Gamma(s)} \ast \chi^{\gamma^G}_{g, g'}) &= (\chi^{\gamma^G}_{(s \rhd_2 g^{-1})^{-1}, ((g^{-1} \rhd_1 s) \rhd_2 g'^{-1})^{-1}} \ast \id_{(g' (s \rhd_2 g^{-1})) \rhd_1 s}) \\
		&\qquad (\id_{\gamma^G((s \rhd_2 g^{-1})^{-1})} \ast \sigma_{g', g^{-1} \rhd_1 s}) (\sigma_{g, s} \ast \id_{\gamma^G(g')}) \\
		\sigma_{g,e} (\iota^{\gamma^\Gamma} \ast \id_{\gamma^G(g)}) &= \id_{\gamma^G(g)} \ast \iota^{\gamma^\Gamma} \\
		\sigma_{e,s} (\id_{\gamma^\Gamma(s)} \ast \iota^{\gamma^G}) &= \iota^{\gamma^G} \ast \id_{\gamma^\Gamma(s)}
	\end{align*}
	for all $g, g' \in G$ and $s, s' \in \Gamma$. Then, $\calc$ becomes a $(G \bowtie \Gamma, G \times \Gamma)$-crossed monoidal category.

	\begin{proof}
		Put $\gamma(g,s) \coloneqq \gamma^G(g) \gamma^\Gamma(s)$ for $(g,s) \in G \bowtie \Gamma$. Put
		\begin{align*}
			\chi^\gamma_{(g, s), (g',s')} \coloneqq (\chi^{\gamma^G}_{g, (s \rhd_2 g'^{-1})^{-1}} \ast \chi^{\gamma^\Gamma}_{g'^{-1} \rhd_1 s,s'}) (\id_{\gamma^G(g)} \ast \sigma^{g', s} \ast \id_{\gamma^\Gamma(s')})
		\end{align*}
		for $(g, s), (g', s') \in G \bowtie \Gamma$. Put $\iota^\gamma \coloneqq \iota^{\gamma^G} \ast \iota^{\gamma^\Gamma}$. Then, $\gamma$ is a $G \bowtie \Gamma$-action on the category $\calc$ as in the proof that the direct product of two algebras is again an algebra. Namely, $\chi^{\gamma^G}$, $\chi^{\gamma^\Gamma}$, $\iota^{\gamma^G}$ and $\iota^{\gamma^\Gamma}$ correspond to the products and units of the algebras. It corresponds to associativity and unit property that $\gamma^G$ and $\gamma^\Gamma$ are group actions. The natural isomorphism $\sigma^{g,s}$ corresponds to the swap between algebras and the last four axioms for $\sigma$ correspond to the commutativity between the swap and the algebra structures. Then, the axioms for $\chi^\gamma$ and $\iota^\gamma$ follow, respectively, as in the proof of the associativity and unit property of the direct product algebra.

		For $\lambda \in \obj(\calc)$, $\mu \in \homog(\calc)$ and $(g,s) \in G \bowtie \Gamma$, put
		\begin{align*}
			J^{\gamma(g,s)}_{\lambda, \mu} &\coloneqq {}^g J^{\gamma^\Gamma(s)}_{\lambda, \mu} J^{\gamma^G(g)}_{{}^{\partial \mu \rhd_2^\Gamma s} \lambda, {}^s \mu} : {}^{h \tilde{\rhd}_2 (g,s)} \lambda {}^{(g,s)} \mu \cong {}^{(g,s)} (\lambda \mu) \\
			\varphi^{\gamma(g,s)} &\coloneqq {}^g \varphi^{\gamma^\Gamma(s)} \varphi^{\gamma^G(g)} : \mathbf{1}_\calc \cong {}^{(g,s)} \mathbf{1}_\calc.
		\end{align*}
		Then, axiom 2 in Definition \ref{definition_crossed_category} follows as in the proof that the composition of monoidal functors is again a monoidal functor. 
		
		We check axiom 3. Let $(\sigma_{\partial \mu \rhd_2^G g,(g \rhd_1^G \partial \mu) \rhd_2^\Gamma s})_\lambda \otimes (\sigma_{g,s})_\mu$ be denoted by a crossing from $s g$ to $(s \rhd_2 g^{-1})^{-1} (g^{-1} \rhd_1 s)$. Then, the first axiom for $\sigma$ is graphically represented as the Yang--Baxter relation in Figure \ref{figure_first_axiom_sigma}. Therefore, the first condition of axiom 3 follows from the graphical calculation in Figure \ref{figure_main1_proof}.
		
		\begin{figure}[htb]
			\centering
			\begin{tikzpicture}
				\draw[->] (2,0) -- (2,-1) -- (0,-3);
				\draw[->,cross] (1,0) -- (0,-1) -- (0,-2) -- (1,-3);
				\draw[->,cross] (0,0) -- (2,-2) -- (2,-3); 
				\node at (0,0.25){$\otimes$};
				\node at (2,-3.25){$\otimes$};
				\node at (1,0.25){$s$};
				\node at (2,0.25){$g$};
				\node at (-0.6,-3.25){$(s \rhd_2 g^{-1})^{-1}$};
				\node at (1.1,-3.25){$g^{-1} \rhd_1 s$};
				\node at (3,-1.5){$=$};
				\begin{scope}[shift={(4,0)}]
					\draw[->] (2,0) -- (0,-2) -- (0,-3);
					\draw[->,cross] (1,0) -- (2,-1) -- (2,-2) -- (1,-3);
					\draw[->,cross] (0,0) -- (0,-1) -- (2,-3); 
					\node at (0,0.25){$\otimes$};
					\node at (2,-3.25){$\otimes$};
					\node at (1,0.25){$s$};
					\node at (2,0.25){$g$};
					\node at (-0.6,-3.25){$(s \rhd_2 g^{-1})^{-1}$};
					\node at (1.1,-3.25){$g^{-1} \rhd_1 s$};
				\end{scope}
			\end{tikzpicture}
			\caption{The first axiom for $\sigma$}
			\label{figure_first_axiom_sigma}
		\end{figure}

		\begin{figure}[htb]
			\centering 
			\begin{tikzpicture}
				\draw (0,0) -- (0,-0.5) arc (180:360:0.25) -- (1,0);
				\draw[cross] (0.5,0) -- (1,-0.5) arc (180:360:0.25) -- (1.5,0);
				\draw[->] (0.25,-0.75) -- (0.25,-1) -- (-1.25,-2.5);
				\draw[->] (1.25,-0.75) -- (1.25,-2) -- (0.75,-2.5);
				\draw[->,cross] (-1.25,0) -- (1.25,-2.5);
				\node at (-1.25,0.25){$\otimes$};
				\node at (1.5,-2.75){$\otimes$};
				\node at (-2,-2.75){$g (s \rhd_2 g'^{-1})^{-1}$};
				\node at (0.25,-2.75){$(g'^{-1} \rhd_1 s)s'$};
				\node at (0,0.25){$g$};
				\node at (0.5,0.25){$s$};
				\node at (1,0.25){$g'$};
				\node at (1.5,0.25){$s'$};
				\node at (2.5,-1.25){$=$};
				\begin{scope}[shift={(5,0)}]
					\draw (0,0) -- (-1.25,-1.25) -- (-1.25,-2) arc (180:360:0.25) -- (-0.75,-1.75) -- (1,0);
					\draw[cross] (0.5,0) -- (1,-0.5) -- (-0.5,-2) arc (180:360:0.25) -- (1.5,-0.5) -- (1.5,0);
					\draw[->] (-1,-2.25) -- (-1,-2.5);
					\draw[->] (-0.25,-2.25) -- (-0.25,-2.5);
					\draw[->,cross] (-1.25,0) -- (1.25,-2.5);
					\node at (-1.25,0.25){$\otimes$};
					\node at (1.5,-2.75){$\otimes$};
					\node at (-2,-2.75){$g (s \rhd_2 g'^{-1})^{-1}$};
					\node at (0.25,-2.75){$(g'^{-1} \rhd_1 s)s'$};
					\node at (0,0.25){$g$};
					\node at (0.5,0.25){$s$};
					\node at (1,0.25){$g'$};
					\node at (1.5,0.25){$s'$};
					\node at (2.5,-1.25){$=$};
				\end{scope}
				\begin{scope}[shift={(10,0)}]
					\draw (0,0) -- (-1.25,-1.25) -- (-1.25,-2) arc (180:360:0.25) -- (-0.75,-1.75) -- (1,0);
					\draw[cross] (0.5,0) -- (-0.75,-1.25) -- (-0.5,-1.5) -- (-0.5,-2) arc (180:360:0.25) -- (1.5,-0.5) -- (1.5,0);
					\draw[->] (-1,-2.25) -- (-1,-2.5);
					\draw[->] (-0.25,-2.25) -- (-0.25,-2.5);
					\draw[->,cross] (-1.25,0) -- (1.25,-2.5);
					\node at (-1.25,0.25){$\otimes$};
					\node at (1.5,-2.75){$\otimes$};
					\node at (-2,-2.75){$g (s \rhd_2 g'^{-1})^{-1}$};
					\node at (0.25,-2.75){$(g'^{-1} \rhd_1 s)s'$};
					\node at (0,0.25){$g$};
					\node at (0.5,0.25){$s$};
					\node at (1,0.25){$g'$};
					\node at (1.5,0.25){$s'$};
				\end{scope}
			\end{tikzpicture}
			\caption{The proof of axiom 3}
			\label{figure_main1_proof}
		\end{figure}

		Next, the second condition of axiom 3 follows since
		\begin{align*}
			&(\chi^\gamma_{(g,s),(g',s')})_{\mathbf{1}_{\calc}} {}^{(g,s)} \varphi^{\gamma(g',s')} \varphi^{\gamma(g,s)} \\
			&= {}^{g(s \rhd_2 g'^{-1})^{-1}} (\chi^{\gamma^\Gamma}_{g'^{-1} \rhd_1 s, s'})_{\mathbf{1}_\calc} (\chi^{\gamma^G}_{g,(s \rhd_2 g'^{-1})^{-1}})_{{}^{g'^{-1} \rhd_1 s}({}^{s'} \mathbf{1}_\calc)} \\
			&\qquad {}^g (\sigma_{g',s})_{{}^{s'} \mathbf{1}_\calc} {}^{g} ({}^s ({}^{g'} \varphi^{\gamma^\Gamma(s')} \varphi^{\gamma^G(g')})) {}^{g} \varphi^{\gamma^\Gamma(s)} \varphi^{\gamma^G(g)} \\
			&= {}^{g(s \rhd_2 g'^{-1})^{-1}} (\chi^{\gamma^\Gamma}_{g'^{-1} \rhd_1 s, s'})_{\mathbf{1}_\calc} (\chi^{\gamma^G}_{g,(s \rhd_2 g'^{-1})^{-1}})_{{}^{g'^{-1} \rhd_1 s}({}^{s'} \mathbf{1}_\calc)} {}^g ({}^{(s \rhd_2 g'^{-1})^{-1}} ({}^{g'^{-1} \rhd_1 s} \varphi^{\gamma^\Gamma(s')})) \\
			&\qquad {}^g ((\sigma_{g',s})_{ \mathbf{1}_\calc} {}^s \varphi^{\gamma^G(g')} \varphi^{\gamma^\Gamma(s)}) \varphi^{\gamma^G(g)} \\
			&= {}^{g(s \rhd_2 g'^{-1})^{-1}}((\chi^{\gamma^\Gamma}_{g'^{-1} \rhd_1 s, s'})_{\mathbf{1}_\calc} {}^{g'^{-1} \rhd_1 s} \varphi^{\gamma^\Gamma(s')})  (\chi^{\gamma^G}_{g,(s \rhd_2 g'^{-1})^{-1}})_{{}^{g'^{-1} \rhd_1 s}\mathbf{1}_\calc} \\
			&\qquad {}^g ({}^{(s \rhd_2 g'^{-1})^{-1}} \varphi^{\gamma^\Gamma(g'^{-1} \rhd_1 s)} \varphi^{\gamma^G((s \rhd_2 g'^{-1})^{-1})} ) \varphi^{\gamma^G(g)} \\
			&= {}^{g(s \rhd_2 g'^{-1})^{-1}}((\chi^{\gamma^\Gamma}_{g'^{-1} \rhd_1 s, s'})_{\mathbf{1}_\calc} {}^{g'^{-1} \rhd_1 s} \varphi^{\gamma^\Gamma(s')} \varphi^{\gamma^\Gamma(g'^{-1} \rhd_1 s)}) \\
			&\qquad (\chi^{\gamma^G}_{g,(s \rhd_2 g'^{-1})^{-1}})_{\mathbf{1}_\calc} {}^g   \varphi^{\gamma^G((s \rhd_2 g'^{-1})^{-1})} \varphi^{\gamma^G(g)} \\
			&= {}^{g(s \rhd_2 g'^{-1})^{-1}} \varphi^{\gamma^\Gamma((g'^{-1} \rhd_1 s)s')} \varphi^{\gamma^G(g (s \rhd_2 g'^{-1})^{-1})} = \varphi^{\gamma(g (s \rhd_2 g'^{-1})^{-1}, (g'^{-1} \rhd_1 s)s')}
		\end{align*}
		by the naturality of $\sigma_{g',s}$ and $\chi^{\gamma^G}_{g,(s \rhd_2 g'^{-1})^{-1}}$, the second axiom for $\sigma$ and the second condition of axiom 3 for $\gamma^G$ and $\gamma^\Gamma$. Then, the third condition of axiom 3 for $\gamma$ follows since
		\begin{align*}
			\iota^{\gamma}_{\lambda \mu} &= \gamma^G(e)(\iota^{\gamma^\Gamma}_{\lambda \mu}) \iota^{\gamma^G}_{\lambda \mu} = \gamma^G(e)(J^{\gamma^\Gamma(e)}_{\lambda, \mu} (\iota^{\gamma^\Gamma}_\lambda \otimes \iota^{\gamma^\Gamma}_{\mu})) J^{\gamma^G(e)}_{\lambda, \mu} (\iota^{\gamma^G}_\lambda \otimes \iota^{\gamma^G}_\mu) \\
			&= \gamma^G(e) J^{\gamma^\Gamma(e)}_{\lambda, \mu} J^{\gamma^G(e)}_{\gamma^\Gamma(e)(\lambda),\gamma^\Gamma(e)(\mu)} (\gamma^G(e)(\iota^{\gamma^\Gamma}_\lambda) \otimes \gamma^G(e)(\iota^{\gamma^\Gamma}_\mu)) = J^{\gamma(e)}_{\lambda, \mu} (\iota^\gamma_\lambda \otimes \iota^\gamma_\mu)
		\end{align*}
		by the third condition of axiom 3 for $\gamma^G$ and $\gamma^\Gamma$ and the naturality of $J^{\gamma^G(e)}$. Finally, the fourth condition of axiom 3 for $\gamma$ follows since $\iota^\gamma_{\mathbf{1}_\calc} = \gamma^G(e)(\iota^{\gamma^\Gamma}_{\mathbf{1}_\calc}) \iota^{\gamma^G}_{\mathbf{1}_\calc} = \gamma^G(e)(\varphi^{\gamma^\Gamma(e)}) \varphi^{\gamma^G(e)} = \varphi^{\gamma(e)}$ by the fourth condition in axiom 3 for $\gamma^G$ and $\gamma^\Gamma$.
	\end{proof}
\end{lemma}

The following theorem is the main result of this article.

\begin{theorem}
	\label{theorem_crossed_center}
	The $(G,\Gamma)$-crossed center $\rmz^{G}_\Gamma(\calc)$ of a non-singular $(G,\Gamma)$-crossed pivotal tensor category $\calc$ becomes a $(G \bowtie \Gamma, G \times \Gamma)$-braided pivotal tensor category. 

	\begin{proof}
		Let $\gamma^G$ (resp. $\gamma^\Gamma$) denote the actions of $G$ (resp. $\Gamma$) with crossed structures given by Lemma \ref{lemma_center_g_action} (resp. Lemma \ref{lemma_center_gamma_action}). For $\lambda \in \obj(\rmz^G_\Gamma(\calc)_{(h,t)})$ and $(g,s) \in G \bowtie \Gamma$, we define $\sigma_{g,s}$ by Figure \ref{figure_definition_sigma}, where we put $d_s \coloneqq d_{\zeta_s}$ for the fixed family $\{ \zeta_s \}_{s \in \Gamma}$ of homogeneous objects with nonzero dimensions.
		
		\begin{figure}[htb]
			\centering 
			\begin{tikzpicture}
				\draw (-0.5,0) arc (180:360:0.5);
				\draw (-0.5,-1.25) arc (180:0:0.5);
				\draw[cross,->] (0,0.75) -- (0,-2.25);
				\draw[fill=white] (-0.5,0) -- (0.5,0) -- (0,0.5) -- cycle;
				\node[block] at(0,-1.5){${}^{(s \rhd_2 g^{-1})^{-1}} r^{g^{-1} \rhd_1 s}_\lambda$};
				\node at(0,1){${}^s ({}^g \lambda)$};
				\node at(0.75,0){$s$};
				\node at(-1,-0.75){$d_{s}^{-1}$};
				\node at(0,-2.5){${}^{(s \rhd_2 g^{-1})^{-1}}({}^{g^{-1} \rhd_1 s} \lambda)$};
			\end{tikzpicture}
			\caption{$(\sigma_{g,s})_\lambda$}
			\label{figure_definition_sigma}
		\end{figure}

		These $\sigma_{g,s}$'s are natural isomorphisms $\gamma^\Gamma(s) \gamma^G(g) \cong \gamma^G((s \rhd_2 g^{-1})^{-1}) \gamma^\Gamma(g^{-1} \rhd_1 s)$ and satisfy the conditions in Lemma \ref{lemma_bicross_action} by the same argument as in the proof that $\gamma^\Gamma$ does not depend on a choice of a family $\{ \zeta_s \}_{s \in \Gamma}$ of homogeneous objects with nonzero dimensions. Thus, $\rmz^G_\Gamma(\calc)$ becomes a $(G,\Gamma)$-crossed pivotal tensor category.

		Then, we endow $\rmz^G_\Gamma(\calc)$ with a $(G,\Gamma)$-braiding. For $\lambda, \mu \in \homog(\rmz^G_\Gamma(\calc))$, we define the isomorphism $b_{\lambda, \mu}$ by Figure \ref{figure_braiding}.
		\begin{figure}[htb]
			\centering
			\begin{tikzpicture}
				\draw (0,0) arc (180:360:0.25);
				\draw (1,0) arc (180:0:0.25);
				\draw[cross,->] (1.25,0.5) -- (1.25,-0.5) -- (0.25,-1.5);
				\draw[cross,->] (0.25,0.5) --(0.25,-0.5) -- (1.25,-1.5);
				\draw[fill=white] (1,0) -- (1.5,0) -- (1.25,-0.25) -- cycle;
				\draw[fill=white] (0,0) -- (0.5,0) -- (0.25,0.25) -- cycle;
				\node at (0.25,0.75){${}^{\partial_\Gamma \mu} \lambda$};
				\node at (1.25,0.75){$\mu$};
				\node at (0.25,-1.75){${}^{\partial_G \lambda} ({}^e \mu)$};
				\node at (1.25,-1.75){$\lambda$};
				\node at (1.75,0){$e$};
				\node at (-0.5,0){$\partial_\Gamma \mu$};
				\node at (-0.5,-0.75){$d_{\partial_\Gamma \mu}^{-1}$};
			\end{tikzpicture}
			\caption{The $(G \bowtie \Gamma, G \times \Gamma)$-braiding on $\rmz^G_\Gamma(\calc)$}
			\label{figure_braiding}
		\end{figure}

		Then, by the Yang--Baxter equation and the definition of $\gamma^\Gamma$, $b_{\lambda, \mu}$ is a morphism of $\rmz^G_\Gamma(\calc)$. Moreover, it is easy to see that the morphism of Figure \ref{figure_braiding_inv} is the inverse of $b_{\lambda, \mu}$.
		\begin{figure}[htb]
			\centering
			\begin{tikzpicture}
				\begin{scope}[yscale=-1]
				\draw (0,0) arc (180:360:0.25);
				\draw (1,0) arc (180:0:0.25);
				\draw[cross,->] (1.25,0.5) -- (1.25,-0.5) -- (0.25,-1.5);
				\draw[cross,->] (0.25,0.5) --(0.25,-0.5) -- (1.25,-1.5);
				\draw[fill=white] (1,0) -- (1.5,0) -- (1.25,-0.25) -- cycle;
				\draw[fill=white] (0,0) -- (0.5,0) -- (0.25,0.25) -- cycle;
				\node at (0.25,0.75){${}^{\partial_\Gamma \mu} \lambda$};
				\node at (1.25,0.75){$\mu$};
				\node at (0.25,-1.75){${}^{\partial_G \lambda} ({}^e \mu)$};
				\node at (1.25,-1.75){$\lambda$};
				\node at (1.75,0){$e$};
				\node at (-0.5,-0.75){$d_{e}^{-1}$};
				\node at (-0.5,0){$\partial_\Gamma \mu$};
				\end{scope}
			\end{tikzpicture}
			\caption{$b^{-1}_{\lambda,\mu}$}
			\label{figure_braiding_inv}
		\end{figure}

		We show the first axiom for $b$. Put $(h,t) \coloneqq \partial \lambda$ and $(k,u) \coloneqq \partial \mu$. For $(g,s) \in G \bowtie \Gamma$, we obtain
		\begin{align*}
			&(\chi^\gamma_{(k,u) \tilde{\rhd}_2 (g,s), \phi(k,u)})_\lambda^{-1} (\chi^\gamma_{\phi((g,s) \tilde{\rhd}_1 (k,u)), (g,s)})_\lambda \\ 
			&=  {}^{((k \rhd_1 s)us^{-1}) \rhd_2 g} ((\chi^{\Gamma}_{k \rhd_1 s , t})^{-1}_\lambda (\chi^\Gamma_{(k \rhd_1 s)us^{-1},s})_\lambda) (\sigma_{g,g \rhd_1 ((k \rhd_1 s)us^{-1})})_{{}^s \lambda}
		\end{align*}
		by the definition of $\chi^\gamma$ and the coherence for the $G$-action $\gamma^G$. Then, the left-hand side of the first axiom for $b$ is equal to 
		\begin{align}
			&d_{(g \rhd_1 s^{-1})^{-1}}^{-2} d_{g \rhd_1 ((k \rhd_1 s)us^{-1})}^{-1} (\sigma_{s^{-1} \rhd_2 g,(g \rhd_1 s^{-1})^{-1}})_{{}^h ({}^e \mu) \lambda} f \nonumber \\
			\label{proof_braiding_first_morphism}
			&((\sigma_{g,g \rhd_1 ((k \rhd_1 s)us^{-1})})_{{}^u \lambda} {}^{g \rhd_1 ((k \rhd_1 s)us^{-1})} (\sigma_{s^{-1} \rhd_2 g,(g \rhd_1 s^{-1})^{-1}})_\lambda \otimes (\sigma_{s^{-1} \rhd_2 g,(g \rhd_1 s^{-1})^{-1}})_\mu {}^g({}^s \iota_\mu^{\Gamma}))
		\end{align}
		for the morphism $f$ defined by Figure \ref{figure_proof_braiding_first}. Here we used the following simple replacement rules that follow from the definition of $\sigma$: when $g \in G$ acts on the triangle representing a retract labeled by $s \in G$, we replace the label with $(s \rhd_2 g^{-1})^{-1}$ and concatenate $(\sigma_{s^{-1} \rhd_2 g,(g \rhd_1 s^{-1})^{-1}})_\lambda$ from below. When $g$ acts on the triangle representing a section labeled by $s$, we replace the label in the same way, concatenate $(\sigma_{s^{-1} \rhd_2 g,(g \rhd_1 s^{-1})^{-1}})_\lambda^{-1}$ from above and multiply the factor $d_s/d_{(g \rhd_1 s^{-1})^{-1}}$.

		\begin{figure}[htb]
			\centering 
			\begin{tikzpicture}
				\draw (-0.25,0) arc (180:360:0.25);
				\draw (-0.5,0.5) -- (-0.5,0) arc (180:360:0.5) -- (0.5,0.5);
				\draw (1,0.5) arc (180:360:0.5);
				\draw (0,-1.5) arc (180:0:0.75);
				\draw[->] (0.75,-2) -- (0.75,-2.25);
				\draw[cross] (1.5,1.25) -- (1.5,-0.5) -- (0.5,-1.5);
				\draw[cross] (0,1.25) -- (0,-0.5) -- (1,-1.5);
				\draw[fill=white] (-0.25,0) -- (0.25,0) -- (0,0.25) -- cycle;
				\draw[fill=white] (-0.5,0.5) -- (0.5,0.5) -- (0,1) -- cycle;
				\draw[fill=white] (1,0.5) -- (2,0.5)-- (1.5,1) -- cycle;
				\draw[fill=white] (0,-1.5) -- (1.5,-1.5) -- (0.75,-2) -- cycle;
				\node at(-2,1.5){${}^{(g \rhd_1 ((k \rhd_1 s)us^{-1}))(g \rhd_1 s^{-1})^{-1}} \lambda$};
				\node at(2,1.5){${}^{(g \rhd_1 s^{-1})^{-1}}({}^{s^{-1} \rhd_2 g}({}^e\mu))$};
				\node at(0.75,-2.5){${}^{(g \rhd_1 s^{-1})^{-1}} ({}^h \mu \lambda)$};
				\node at(3,0.5){$(g \rhd_1 s^{-1})^{-1}$};
				\node at(2.5,-1.5){$(g \rhd_1 s^{-1})^{-1}$};
			\end{tikzpicture}
			\caption{The first axiom for $b$}
			\label{figure_proof_braiding_first}
		\end{figure}

		On the other hand, we have
		\begin{align*}
			&(\chi^{\gamma}_{(h,t) \tilde{\rhd}_2 (g,s),\psi(h,t)})^{-1}_{\mu} (\chi^\gamma_{\psi((g,s) \tilde{\rhd}_1 (h,t)),(g,s)})_{\mu} \\
			&= {}^{((h \rhd_1 s)ts^{-1}) \rhd_2 g} (\sigma_{h,h \rhd_1 s})_{{}^e \mu}^{-1} {}^{(((h \rhd_1 s)ts^{-1}) \rhd_2 g)(s \rhd_2 h)}({}^s \iota^\Gamma_\mu {}^{g^{-1}} \iota^{\Gamma-1}_{{}^g({}^s \mu)})
		\end{align*}
		by the definition of $\chi^\gamma$ and the coherence. Then, by the naturality of $b$ and the replacement rules, we see that the right-hand side of the first axiom for $b$ is equal to the morphism (\ref{proof_braiding_first_morphism}).

		Then, we show the second axiom for $b$. The left-hand side of the axiom is equal to the morphism in Figure \ref{figure_proof_braiding_second} up to the factor $d_{\partial_G \lambda_2 \rhd_1 \partial_\Gamma \mu}^{-1} d_{\partial_\Gamma \mu}^{-1}$ for $\lambda_1, \lambda_2, \mu \in \homog(\rmz_\Gamma^G(\calc))$ by the definition of $J^{\gamma}$. On the other hand, by the naturality of $h^{\lambda_1}$, the right-hand side of the axiom is equal to the composition of $d_{\partial_\Gamma \mu}^{-2} \chi^\gamma_{(\partial_G \lambda_1 ,e),(\partial_G \lambda_2, e)} {}^{\partial_G \lambda_2} \iota^\Gamma_{{}^{\partial_G \lambda_2}({}^e \mu)}$ and the morphism in Figure \ref{figure_proof_braiding_second}. Then, the second axiom follows since $\chi^\gamma_{(\partial_G \lambda_1 ,e),(\partial_G \lambda_2, e)} = {}^{\partial_G \lambda_1} \iota^{\Gamma-1}_{{}^{\partial_G \lambda_2} ({}^e \mu)}$ by the definition of $\chi^\gamma$ and the coherence for the $\Gamma$-action $\gamma^\Gamma$.
		\begin{figure}[htb]
			\centering
			\begin{tikzpicture}
				\draw (0,0) arc (180:360:0.25);
				\draw (1,0) arc (180:0:0.25);
				\draw (-1,0) arc (180:360:0.25);
				\draw[cross,->] (1.25,0.5) -- (1.25,-0.5) -- (-0.75,-2.5);
				\draw[cross,->] (0.25,0.5) --(0.25,-0.5) -- (1.25,-1.5) -- (1.25,-2.5);
				\draw[cross,->] (-0.75,0.5) -- (-0.75,-1.5) --(0.25,-2.5);
				\draw[fill=white] (1,0) -- (1.5,0) -- (1.25,-0.25) -- cycle;
				\draw[fill=white] (0,0) -- (0.5,0) -- (0.25,0.25) -- cycle;
				\draw[fill=white] (-1,0) -- (-0.5,0) -- (-0.75,0.25) -- cycle;
				\node at (0.25,0.75){${}^{\partial_\Gamma \mu} \lambda_2$};
				\node at (-1.25,0.75){${}^{\partial \lambda_2 \tilde{\rhd}_2 \partial_\Gamma \mu} \lambda_1$};
				\node at (1.25,0.75){$\mu$};
				\node at (0.25,-2.75){$\lambda_1$};
				\node at (1.25,-2.75){$\lambda_2$};
				\node at (-1.5,-2.75){${}^{\partial_G \lambda_1 \lambda_2}({}^e \mu)$};
			\end{tikzpicture}
			\caption{The second axiom for $b$}
			\label{figure_proof_braiding_second}
		\end{figure}

		Finally, we show the third axiom for $b$. By definition, the left-hand side of the third axiom is equal to the leftmost side in Figure \ref{figure_proof_braiding_third} up to the factor $d_{\partial_\Gamma \mu_1}^{-1} d_{\partial_\Gamma \mu_2}^{-1}$. Since the right-hand side of the third axiom is equal to the rightmost side in Figure \ref{figure_proof_braiding_third} up to the same factor,
		 the third axiom follows from the calculation in Figure \ref{figure_proof_braiding_third}.
		\begin{figure}[htb]
			\centering 
			\begin{tikzpicture}
				\draw (-0.25,0) arc (180:360:0.25);
				\draw (-0.5,0.5) -- (-0.5,0) arc (180:360:0.5) -- (0.5,0.5);
				\draw[->] (1.5,-0.5) -- (0,-2);
				\draw (1,0) arc (180:0:0.5);
				\draw[cross, ->] (0,1.25) -- (0,-1) -- (1,-2);
				\draw[cross] (1.25,1.25) -- (1.25,0);
				\draw[cross] (1.75,1.25) -- (1.75,0);
				\draw[fill=white] (-0.25,0) -- (0.25,0) -- (0,0.25) -- cycle;
				\draw[fill=white] (-0.5,0.5) -- (0.5,0.5) -- (0,1) -- cycle;
				\draw[fill=white] (1,0) -- (2,0) -- (1.5,-0.5) -- cycle;
				\node at(0,1.5){${}^{\partial_\Gamma \mu_1 \mu_2} \lambda$};
				\node at(1.25,1.5){$\mu_1$};
				\node at(1.75,1.5){$\mu_2$};
				\node at(-0.25,-2.25){${}^{\partial_G \lambda} ({}^e(\mu_1 \mu_2))$};
				\node at(1,-2.25){$\lambda$};
				\node at(-1,0.5){$\partial_\Gamma \mu_1$};
				\node at(2.25,0){$e$};
				\node at(2.5,-0.25){$=$};
				\begin{scope}[shift={(4,0)}]
					\draw (-0.25,0) arc (180:360:0.25);
					\draw (-0.5,0.5) -- (-0.5,0) arc (180:360:0.5) -- (0.5,0.5);
					\draw[->] (1.5,-0.5) -- (0,-2);
					\draw (1,0) arc (180:0:0.5);
					\draw (1,1) arc (180:360:0.25);
					\draw (1.5,1) arc (180:360:0.25);
					\draw (1,1.5) arc (180:0:0.25);
					\draw (1.5,1.5) arc (180:0:0.25);
					\draw[cross, ->] (0,2) -- (0,-1) -- (1,-2);
					\draw[cross] (1.25,2) -- (1.25,0);
					\draw[cross] (1.75,2) -- (1.75,0);
					\draw[fill=white] (-0.25,0) -- (0.25,0) -- (0,0.25) -- cycle;
					\draw[fill=white] (-0.5,0.5) -- (0.5,0.5) -- (0,1) -- cycle;
					\draw[fill=white] (1,0) -- (2,0) -- (1.5,-0.5) -- cycle;
					\draw[fill=white] (1,1) -- (1.5,1) -- (1.25,1.25) -- cycle;
					\draw[fill=white] (1.5,1) -- (2,1) -- (1.75,1.25) -- cycle;
					\draw[fill=white] (1,1.5) -- (1.5,1.5) -- (1.25,1.25) -- cycle;
					\draw[fill=white] (1.5,1.5) -- (2,1.5) -- (1.75,1.25) -- cycle;
					\node at(0,2.25){${}^{\partial_\Gamma \mu_1 \mu_2} \lambda$};
					\node at(1.25,2.25){$\mu_1$};
					\node at(1.75,2.25){$\mu_2$};
					\node at(-0.25,-2.25){${}^{\partial_G \lambda} ({}^e(\mu_1 \mu_2))$};
					\node at(1,-2.25){$\lambda$};
					\node at(-1,0.5){$\partial_\Gamma \mu_1$};
					\node at(2.25,0){$e$};
					\node at(2.25,1.5){$e$};
					\node at(-1,-0.25){$d_e^{-2}$};
				\end{scope}
				\begin{scope}[shift={(0,-5)}]
					\draw (-0.25,0) arc (180:360:0.25);
					\draw (-0.5,0.5) -- (-0.5,0) arc (180:360:0.5) -- (0.5,0.5);
					\draw[->] (0,-1.5) -- (0,-2);
					\draw (1,1.5) arc (180:0:0.25);
					\draw (1.5,1.5) arc (180:0:0.25);
					\draw[cross] (1.25,2) -- (1.25,0) -- (-0.25,-1.5);
					\draw[cross] (1.75,2) -- (1.75,0) -- (0.25,-1.5);
					\draw[cross, ->] (0,2) -- (0,0) -- (2,-2);
					\draw[fill=white] (-0.25,0) -- (0.25,0) -- (0,0.25) -- cycle;
					\draw[fill=white] (-0.5,0.5) -- (0.5,0.5) -- (0,1) -- cycle;
					\draw[fill=white] (1,1.5) -- (1.5,1.5) -- (1.25,1.25) -- cycle;
					\draw[fill=white] (1.5,1.5) -- (2,1.5) -- (1.75,1.25) -- cycle;
					\node at(0,2.25){${}^{\partial_\Gamma \mu_1 \mu_2} \lambda$};
					\node at(1.25,2.25){$\mu_1$};
					\node at(1.75,2.25){$\mu_2$};
					\node at(-0.25,-2.25){${}^{\partial_G \lambda} ({}^e(\mu_1 \mu_2))$};
					\node at(2,-2.25){$\lambda$};
					\node at(-1,0.5){$\partial_\Gamma \mu_1$};
					\node at(2.25,1.5){$e$};
					\node[block] at(0,-1.4){${}^{\partial_G \lambda} J^{\gamma^\Gamma(e)}_{\mu_1,\mu_2}$};
					\node at(-1,-0.25){$=$};
				\end{scope}
				\begin{scope}[shift={(4,-5)}]
					\draw (0.75,0) arc (180:360:0.25);
					\draw (-0.5,1.25) arc (180:360:0.5);
					\draw[->] (0,-1.5) -- (0,-2);
					\draw (1,1.5) arc (180:0:0.25);
					\draw (1.5,1.5) arc (180:0:0.25);
					\draw[cross] (1.25,2) -- (1.25,1) -- (-0.25,-0.5) -- (-0.25,-1.5);
					\draw[cross] (1.75,2) -- (1.75,0) -- (0.25,-1.5);
					\draw[cross, ->] (0,2) -- (0,1) -- (2,-1) -- (2,-2);
					\draw[fill=white] (0.75,0) -- (1.25,0) -- (1,0.25) -- cycle;
					\draw[fill=white] (-0.5,1.25) -- (0.5,1.25) -- (0,1.75) -- cycle;
					\draw[fill=white] (1,1.5) -- (1.5,1.5) -- (1.25,1.25) -- cycle;
					\draw[fill=white] (1.5,1.5) -- (2,1.5) -- (1.75,1.25) -- cycle;
					\node at(0,2.25){${}^{\partial_\Gamma \mu_1 \mu_2} \lambda$};
					\node at(1.25,2.25){$\mu_1$};
					\node at(1.75,2.25){$\mu_2$};
					\node at(-0.25,-2.25){${}^{\partial_G \lambda} ({}^e(\mu_1 \mu_2))$};
					\node at(2,-2.25){$\lambda$};
					\node at(-1,1.25){$\partial_\Gamma \mu_1$};
					\node at(2.25,1.5){$e$};
					\node[block] at(0,-1.4){${}^{\partial_G \lambda} J^{\gamma^\Gamma(e)}_{\mu_1,\mu_2}$};
					\node at(-1.5,-0.25){$=$};
				\end{scope}
			\end{tikzpicture}
			\caption{The third axiom for $b$}
			\label{figure_proof_braiding_third}
		\end{figure}
	\end{proof}
\end{theorem}

\section*{Acknowledgements}

The author thanks his supervisor Yasuyuki Kawahigashi for his constant support and advice. He thanks Yuki Arano, Kan Kitamura, and Toshihiko Masuda for useful comments. This work was supported by RIKEN Junior Research Associate program and JSPS KAKENHI Grant Number JP23KJ0540. He thanks their financial support.

\section*{Declarations}

\textbf{Funding}: This work was supported by RIKEN Junior Research Associate program and JSPS KAKENHI Grant Number JP23KJ0540.

\noindent
\textbf{Conflicts of interest/Competing interests}: The author has no conflicts of interest to declare that are relevant to the content of this article.

\bibliography{oikawa_center}
\bibliographystyle{alpha}

\end{document}